\documentclass{cjour}
\usepackage{graphicx}
\usepackage{amsmath}
\usepackage{amssymb}
\usepackage{latexsym}
\usepackage{url}
\usepackage{pifont}
\usepackage{graphicx}
\usepackage{array}
\newcommand{\nc}{\newcommand}
\nc{\EQ}[1]{(\ref{eq:#1})}
\nc{\TAB}[1]{\ref{tab:#1}}
\nc{\FIG}[1]{\ref{fig:#1}}
\nc{\SEC}[1]{\ref{sec:#1}}
\nc{\mtrx}[1]{\mathsf{#1}}
\nc{\vctr}[1]{\mathsf{#1}}
\nc{\intd}{\mathrm{d}}
\nc{\rmd}{\mathrm{d}}
\nc{\IR}{\mathbb{R}}
\nc{\IN}{\mbox{$\mathbb{N}$}}
\nc{\II}{\mbox{$\mathbb{I}$}}
\nc{\IZ}{\mbox{$\mathbb{Z}$}}
\nc{\IC}{\mathbb{C}}
\nc{\inpr}{\hspace{-1pt}\cdot\hspace{-1pt}}
\nc{\half}{\scriptstyle{\frac{1}{2}}}
\nc{\shalf}{\scriptstyle{\half}} 
\nc{\sign}{\mathrm{sign}}
\nc{\rank}{\mathrm{rank}}
\nc{\diag}{\mathrm{diag}}
\nc{\order}{\mathcal{O}}
\nc{\Ord}{O}
\nc{\ord}{o}

\begin{document}

%%%%% To be entered at Academic Press: =====>>

% \journame{}
% \articlenumber{}
% \yearofpublication{}
% \volume{}
% \cccline{}
% \received{}
% \revised{}
% \accepted{}

\authorrunninghead{VAN BRUMMELEN AND VENNER}
\titlerunninghead{MULTILEVEL EVALUATION OF MULTIDIMENSIONAL INTEGRAL TRANSFORMS}

% communication line, use: \commline{Communicated by...}
% \commline{ }

%\setcounter{page}{261} %% This command is optional.

%% <<== End of commands to be entered at Academic Press

%%  Authors, start here ==>>

%\draft % Optional, will cause a line at the bottom of each page
%% with the words `Draft' and the time and date that the article
%% was LaTeXed. Will also double space text.

\title{Multilevel Evaluation of Multidimensional Integral Transforms with
 Asymptotically Smooth Kernels}

\author{
E.H.\ van Brummelen
}
\affil{
Eindhoven University of Technology, Department of Mechanical Engineering,
P.O. Box 513, 5600 MB Eindhoven, The Netherlands.}

\and 
\author{C.H.\ Venner}
\affil{University of Twente, Faculty of Engineering Technology,
P.O. Box 217, 7500 AE Enschede, The Netherlands.}
\email{e.h.v.brummelen@tue.nl; c.h.venner@utwente.nl.}

%-----------------------[ABSTRACT]-----------------------------
\abstract{
In many practical applications of numerical methods a substantial increase in 
efficiency can be obtained by using local grid refinement, since the solution is 
generally smooth in large parts of the domain and large gradients occur only locally. 
Fast evaluation of integral transforms on such an adaptive grid requires an algorithm 
that relies on the smoothness of the continuum kernel only, independent of its 
discrete form. A multilevel algorithm with this property was presented
in~\cite[A.~Brandt and C.H.~Venner, SIAM J. Sci. Stat. Comput. {\bf 19} (1998) pp.468-492]{Brandt:1998fk}.
Ref.~\cite{Brandt:1998fk} shows that already on a uniform grid the new algorithm is more efficient 
than earlier fast evaluation algorithms, and elaborates the application to
one-dimensional transforms.
The present work analyses the extension and implementation of the algorithm for 
multidimensional transforms. The analysis conveys that the multidimensional extension is 
nontrivial, on account of the occurence of nonlocal corrections. However, by virtue
of the asymptotic smoothness properties of the continuum kernel, these corrections
can again be evaluated fast. By recursion, it is then possible to obtain the optimal
work estimates indicated in~\cite{Brandt:1998fk}.
Currently, only uniform grids are considered. Detailed numerical results will 
be presented for a two dimensional model problem. The results demonstrate that with the new 
algorithm the evaluation of multidimensional transforms is also more efficient 
than with previous algorithms.}
%-----------------------[KEYWORDS ETC.]-------------------------

%-----------------------[INTRO]---------------------------------
\section{Introduction}
\label{c:intro}
In many fields in mathematics, physics and engineering, the numerical
evaluation of integral transforms or multi-integrals of the type:
\begin{equation}\label{eq:int}
Gu(\vctr{x}) \!=\! 
\int_{\Omega}G(\vctr{x},\vctr{y})u(\vctr{y})\,{\mathrm{d}}\vctr{y},
\qquad
\vctr{x}\in\underline{\Omega}\subset\IR^{\underline{d}},
\vctr{y}\in\Omega\subset\IR^{d},
\end{equation}
is a frequently arising task, \mbox{e.g.}, in elasticity-problems, 
integro-differential equations, integral equations, astrophysics and
computer graphics. The evaluation of~\EQ{int} can be a task by
itself or a subtask in the solution of a (system of) integro-differential
equation(s). In the latter case,~$u(\vctr{y})$ is the unknown function.

To evaluate~\EQ{int} numerically, the continuous transform is replaced by a 
matrix multiplication or ``multisummation'', {i.e.} at the expense of a 
discretization error the evaluation of \EQ{int} is replaced by the~$n$-vector
$\vctr{Gu}=\mtrx{G}\inpr\vctr{u}$, given
the~$\bar{n}\times n$ dense matrix~$\mtrx{G}$ and the~$\bar{n}$-vector 
$\vctr{u}$. Multisummations of this form also appear in, for instance,
particle physics (Coulombic molecular interaction).

Straightforward evaluation of the matrix-vector product 
$\mtrx{G}\inpr\vctr{u}$ involves 
$\bar{n}n$ operations. If the matrix~$\mtrx{G}$ has arbitrary entries, no 
faster method than straightforward multiplication exists. However, many cases 
arise in which the ``discrete kernel''~$\mtrx{G}$ has special properties 
that can be used to obtain a fast evaluation algorithm. Several approaches
have been suggested to reduce the computational cost of the multisummation to
below~$\bar{n}n$ operations, by exploiting such special properties, 
\mbox{e.g.}, hierarchical solvers for many body interaction 
problems~\cite{Appel:1985hl,Barnes:1986fp}, multipole expansions~\cite{Nowak:1986nx}, Fast 
Fourier Transform based schemes~\cite{Reichel:1986oq} and wavelet 
techniques. 
%Generally, these techniques have drawbacks, such as a limited 
%accuracy, limitation to potential type kernels, or they require a significant 
%amount of matrix manipulations to arrive at the sparse matrix which enables 
%the fast evaluation. 

In~\cite{Brandt:1990kl}, a general approach referred to as {\em multilevel matrix multiplication} 
or {\em multi\-level multi-integration} was presented. 
The algorithm has been applied to, for instance, 
integral transforms in elastohydrodynamic lubrication problems~\cite{Venner:2000tg} and,
more recently, to integral transforms emanating from discretizations
of the Laplace and Helmholtz equations by the boundary-element method~\cite{Grigoriev:2004hc,Grigoriev:2004jt}.
The algorithm in~\cite{Brandt:1990kl} relies on the smoothness of the matrix~$\mtrx{G}$. For 
particle problems ${G}_{ij}=G({x}_i,{y}_j)$ and the 
smoothness of the discrete kernel follows immediately from the
smoothness of the continuum kernel.
However, the discretization of~\EQ{int} presented in~\cite{Brandt:1990kl}
yields a matrix of which the smoothness is not only determined by 
the smoothness of the continuum kernel, but also by the applied grid. 
The fast evaluation algorithm then requires grid uniformity.

On the other hand, in practical 
applications, {e.g.}, in contact mechanics and in lubrication, a substantial 
increase in efficiency can be obtained by employing non-uniform
grids, since the solution is often smooth in large parts of the
domain and large gradients occur only locally. Moreover, if~$u(\vctr{y})$
has some singularity, local grid refinement is even 
imperative to maintain an efficient work to accuracy relationship.  

The multilevel methodology in principle allows local grid refinements 
in a very natural way, see~\cite{Bai:1986kk,Brandt:1977ad,Brandt:1982bf}, but to 
implement these techniques for integral transforms a new algorithm had to 
be developed. This new algorithm was presented in~\cite{Brandt:1998fk}. For its 
efficiency, the algorithm relies exclusively on the smoothness of the 
{\em continuum} kernel, thereby allowing the use of local grid refinements and
grid adaptivity. In~\cite{Brandt:1998fk}, it was tested for a one dimensional problem 
on a uniform grid and it was shown that already on a uniform grid the 
evaluation is more efficient than with previous algorithms. 
The application to an actual one dimensional problem where local refinement 
is essential to maintain optimal efficiency was discussed in~\cite{Brandt:1998zm}. 

In the present work, the extension and implementation of the algorithm for 
higher dimensional transforms is discussed. To separate the complications
of grid non-uniformity and multiple dimensions, only uniform grids are
considered. The implementation with locally refined grids is deferred to 
future research. Numerical results are presented for a two dimensional model 
problem.

%-----------------------[DISCR]-----------------------------
\section{Discretization}
\label{s:discr}
In this section we briefly review the discretization procedure 
for~\EQ{int}. The details of the procedure can be found in~\cite{Brandt:1998fk}.
Throughout, it will be assumed that~$\underline{d}=d$. The generalization
to more general cases is straightforward. 

The domain~$\Omega$ is divided into subdomains 
$\Omega_{\vctr{j}}^{\vctr{h}}=
\{\vctr{y}\in\IR^d\mid{{y}_k}_{{j}_k}
\leq {y}_k \leq {{y}_k}_{{j}_k+1}, 
1 \leq k \leq d\}$. 
The resulting grid,~$\{\vctr{y}_{\vctr{j}}\}$, is referred to as the integration grid. 
The integral~\EQ{int} can now be rewritten as a 
summation of the contributions of the individual subdomains, defined by: 
\begin{equation}
\label{eq:cont}
G^{\vctr{h}}_{\vctr{j}}u(\vctr{x}) 
=\int_{\Omega_{\vctr{j}}^{\vctr{h}}}
G(\vctr{x},\vctr{y})u(\vctr{y})\,\mathrm{d}\vctr{y}.
\end{equation}
Next, let~$G^{\vctr{l}}(\vctr{x},\vctr{y})$ be a family of kernels, defined recursively, 
\begin{equation}\label{eq:kernfam}
\begin{array}{rl}
G^{{0}}(\vctr{x},\vctr{y})&=G(\vctr{x},\vctr{y}),\\[4mm]
G^{\vctr{l}}(\vctr{x},\vctr{y})&=
\displaystyle{\int_{{x}_k}^{{y}_k}} G^{\vctr{l}-\vctr{e}_k}(\vctr{x},\vctr{y}+(\eta-{y}_k)\,\vctr{e}_k)
\,\mathrm{d}\eta,
\end{array}
\end{equation}
where~$\vctr{e}_k$ denotes the~$k-${th} unit vector. 
Note that~$\vctr{y}-{y}_k\vctr{e}_k$ is~$\vctr{y}$ with its~$k$-th 
component set to zero.
In many practical cases, 
\mbox{e.g.}, for the logarithmic kernel in~\cite{Brandt:1998fk} and for the kernel
in our model problem,
\begin{equation}
\label{eq:modker}
G(\vctr{x},\vctr{y})=|\vctr{y}-\vctr{x}|^{-1},
\qquad{}\vctr{x},\vctr{y}\in\IR^2,
\end{equation}
it is possible to derive~$G^{\vctr{l}}(\vctr{x},\vctr{y})$ analytically. 

The function~$u(\vctr{y})$ is approximated on~$\Omega_{\vctr{j}}^{\vctr{h}}$ by
$\tilde{u}_{\vctr{j}}^{\vctr{h}}(\vctr{y})$, an order 
$\vctr{s}=(s_1,\ldots,s_d)$ interpolation polynomial, i.e. a 
polynomial of degree~$\vctr{s}-1$. The interpolation is done from 
a data-grid of points,~$\{\vctr{z}_{\vctr{j}}\}$, on which for every site
$u^{\vctr{h}}_{\vctr{j}}=u(\vctr{z}^{\vctr{h}}_{\vctr{j}})$ is given. For 
smallest errors, 
the integration interval should be central relative to the interpolation 
points. However, near external boundaries this may no longer be possible. 

A discrete approximation to~\EQ{cont} is obtained by replacing~$u(\vctr{y})$ by
$\tilde{u}^{\vctr{h}}(\vctr{y})$ and integrating by parts~$\vctr{s}$ times:
\begin{equation}\label{eq:pint}
G^{\vctr{h}}_{\vctr{j}}\tilde{u}^{\vctr{h}}(\vctr{x})= 
\sum_{\vctr{l}={1}}^{\vctr{s}}
\sum_{\vctr{a}={0}}^{{1}} (-1)^{|\vctr{l}|+|\vctr{a}|}
G^{\vctr{l}}(\vctr{x},\vctr{y}_{\vctr{j}+\vctr{a}})
\tilde{u}_{\vctr{j}}^{\vctr{h},(\vctr{l}-{1})}
(\vctr{y}_{\vctr{j}+\vctr{a}}),
\end{equation}
where~$|\vctr{l}|=\sum_{k=1}^{d} {l}_k$,~$|\vctr{a}|=\sum_{k=1}^{d} {a}_k$ and 
$\tilde{u}_{\vctr{j}}^{\vctr{h},(\vctr{l}-{1})}(\vctr{y})$ denotes the 
${l}_k-1$ 
derivative of~$\tilde{u}_{\vctr{j}}^{\vctr{h}}(\vctr{y})$ to~$\vctr{y}_k$ for 
all~$k$.  
Note that summation 
over a vector implies summation over each of the components of the vector, 
so that the summation in~\EQ{pint} actually extends over all vertices
of the subdomain~$\Omega_{\vctr{j}}^{\vctr{h}}$. 

The integral transform~\EQ{int} can now be approximated by taking 
the sum of~\EQ{pint} over all subdomains~$\Omega_{\vctr{j}}^{\vctr{h}}$. Rewriting 
this summation, we obtain a sum of~$\prod_{k=1}^{d}s_k$ discrete 
subtransforms,~$S^{\vctr{h},\vctr{l}}(\vctr{x})$, and boundary terms,~$B^{\vctr{h},\vctr{l}}(\vctr{x})$:
\begin{equation}\label{eq:whk}
G^{\vctr{h}}\tilde{u}^{\vctr{h}}(x)= 
\sum_{\vctr{l}={1}}^{\vctr{s}}(-1)^{d+|\vctr{l}|}
B^{\vctr{h},\vctr{l}}(\vctr{x})+
\sum_{\vctr{l}={1}}^{\vctr{s}}(-1)^{d+|\vctr{l}|}
S^{\vctr{h},\vctr{l}}(\vctr{x}),
\end{equation}
with the discrete subtransforms~$S^{\vctr{h},\vctr{l}}(\vctr{x})$ defined by:
\begin{equation}\label{eq:S}
S^{\vctr{h},\vctr{l}}(\vctr{x})=
\sum_{\vctr{j}={0}}^{\vctr{n}}
G^{\vctr{l}}(\vctr{x},\vctr{y}_{\vctr{j}})
U_{\vctr{j}}^{\vctr{h},\vctr{l}},
\end{equation}
where
\begin{equation}\label{eq:U}
U_{\vctr{j}}^{\vctr{h},\vctr{l}}=\left\{
\begin{array}{ll}
\displaystyle{\sum_{\vctr{a}={0}}^{{1}}}(-1)^{d+|\vctr{a}|}\,
\tilde{u}_{\vctr{j}-\vctr{a}}^{\vctr{h},(\vctr{l}-{1})}
(\vctr{y}_{\vctr{j}}),
&\quad
\forall k \,(\frac{1}{2}s_k\leq{j}_k\leq{}n_k-\frac{1}{2}s_k),
\\
0,&\quad \mbox{otherwise}.
\end{array}\right.
\end{equation}
The boundary terms,~$B^{\vctr{h},\vctr{l}}(\vctr{x})$, extend over all nodes bounding the 
domain~$\Omega$ and the subdomains~$\Omega_{\vctr{j}}^{\vctr{h}}$ where the 
integration interval is not central relative to the interpolation points:
\begin{equation}\label{eq:B}
B^{\vctr{h},\vctr{l}}(\vctr{x})=
\sum_{\vctr{j}={0}}^{\vctr{n}}
G^{\vctr{l}}(\vctr{x},\vctr{y}_{\vctr{j}}) 
V_{\vctr{j}}^{\vctr{h},\vctr{l}},
\end{equation}
where
\begin{displaymath}
V_{\vctr{j}}^{\vctr{h},\vctr{l}}=\left\{
\begin{array}{ll}
\displaystyle{\sum_{\vctr{a}={0}}^{{1}}}(-1)^{d+|\vctr{a}|}
\tilde{v}_{\vctr{j}}^{\vctr{h},(\vctr{l}-{1})}(\vctr{y}_{\vctr{j}-\vctr{a}}),
&\quad
\exists k\,(0\leq{j_k}<\frac{1}{2}s_k\,\vee
\,n_k-\frac{1}{2}s_k<{j_k}\leq{n_k}),
\\
0,&\quad
\mbox{otherwise}.
\end{array}\right.
\end{displaymath}
with
\begin{displaymath}
\tilde{v}_{\vctr{j}}^{\vctr{h},(\vctr{l}-{1})}(\vctr{y}_{\vctr{j}}) =\left\{
\begin{array}{ll}
\tilde{u}_{\vctr{j}}^{\vctr{h},(\vctr{l}-{1})}(\vctr{y}_{\vctr{j}}),
& 
\quad\forall k \, (0\leq{j_k}<{n_k}),
\\
0,&\quad\mbox{otherwise}.
\end{array}\right.
\end{displaymath}
Notice that by~\EQ{S}, the kernel in each of the transforms
follows from the continuum kernel~$G(\vctr{x},\vctr{y})$
by integration. 

Assuming that~$u(\vctr{y})$ is~$\vctr{s}$ times 
differentiable on~$\Omega$, in the case of a uniform grid, 
the discretization error, {i.e.} the difference between~\EQ{whk} 
and~\EQ{int}, per unit of integration is bounded by 
\begin{equation}
\label{eq:tau^h}
|G^{\vctr{h}}\tilde{u}^{\vctr{h}}(\vctr{x})-G{u}(\vctr{x})|
%\equiv \tau^{\vctr{h}}(\vctr{x})
\leq \alpha_1\,\sum_{k=1}^{d}\,
(\gamma_1 h_k)^{s_k} 
\big\|u^{(s_k \vctr{e}_k)}\big\|_{\text{max},\Omega}\,
\| G\|_{1,\Omega},
\end{equation}
with~$h_k$ the mesh size of~$\{\vctr{y}_\vctr{j}\}$ in the~$k$-direction,
$\|u^{(s_k \vctr{e}_k)}\|_{\text{max},\Omega}$ the maximum of the 
$s_k$ derivative of~$u(\vctr{y})$ to~$y_k$ on~$\Omega$ and 
$\|G\|_{1,\Omega}$ the average of~$|G(\vctr{x},\vctr{y})|$ over the 
integration 
domain for a particular~$\vctr{x}$. 
Further,~$\alpha_1$ and~$\gamma_1$ are constants, 
the latter depending on the interpolation geometry.

As a result of symmetry and anti-symmetry of the interpolation polynomials, 
on a uniform grid many of the derivatives of the approximating functions,
$\tilde{u}_{\vctr{j}}^{\vctr{h}}(\vctr{y})$, are continuous across the 
subdomain boundaries and, consequently, the 
corresponding~$U_{\vctr{j}}^{\vctr{h},\vctr{l}}$ vanish for all~$\vctr{j}$. 
In particular,~$U_{\vctr{j}}^{\vctr{h},\vctr{l}}=0$ for any~${l}_k$ odd and~$s_k$ even, 
i.e. if the integration grid coincides with the data-grid
(\/${y_k}_{j_k}={z_k}_{j_k}$\/), and for any~$l_k$ even and~$s_k$ 
odd, if the integration grid coincides with data-grid midpoints
(\/${y_k}_{j_k}=({z_k}_{j_k}+{z_k}_{j_k-1})/2$\/); see 
\cite{Brandt:1998fk}. Hence, the number of transforms that actually need to be 
evaluated is~$\prod_{k=1}^{d}\!(\bar{s}_k/2)$, 
where~$\bar{s}_k=s_k$ if 
$s_k$ is even and~$\bar{s}_k=s_k+1$ if~$s_k$ is odd.

To illustrate the discretization procedure, consider the 
discretization~\EQ{whk} with kernel~\EQ{modker} and
$\tilde{u}_\vctr{j}^\vctr{h}(\vctr{y})$ a bi-linear interpolation 
from the data-grid (\/$\vctr{s}=2$\/).
The integration grid and data-grid coincide, as is usual for~$\vctr{s}$ 
is even.
It is easily verified that indeed  
$U^{\vctr{h},\vctr{l}}_\vctr{j}$ vanishes if either~$l_1=1$ or~$l_2=1$. Hence, 
only the discrete transform~$S^{\vctr{h},(2,2)}(\vctr{x})$ needs to be evaluated.
Introducing~$\vctr{t}=\vctr{y}-\vctr{x}$, the integrated kernel for the transform reads:
\begin{multline}
\label{eq:G22}
G^{(2,2)}(\vctr{x},\vctr{y})=
\displaystyle{\frac{1}{2}} {t}_{1} |{t}_{1}| {t}_{2}\:
\mbox{arcsinh}\left(\displaystyle{\frac{{t}_2}{{t}_1}}\right)+
\displaystyle{\frac{1}{2}} {t}_{2} |{t}_{2}| {t}_{1}\:
\mbox{arcsinh}\left(\displaystyle{\frac{{t}_1}{{t}_2}}\right)\\
+\displaystyle{\frac{1}{6}} \left(|{t}_1|^3+|{t}_2|^3 
-\left({t}_1^2+{t}_2^2 \right)^{3/2}\right).
\end{multline}
Inserting the bi-linear approximation into~\EQ{S} and~\EQ{U} yields:
\begin{equation}
\label{eq:S22}
S^{\vctr{h},(2,2)}(\vctr{x})=
\sum_{\vctr{j}={1}}^{\vctr{n}-{1}}
G^{(2,2)}(\vctr{x},\vctr{y}_{\vctr{j}}) U_{\vctr{j}}^{\vctr{h},(2,2)},
\end{equation}
with the stencil of~$U_{\vctr{j}}^{\vctr{h},(2,2)}$ 
(for stencil notation see, e.g.,~\cite{Wesseling:1992rp})
\begin{equation}
U_{\vctr{j}}^{\vctr{h},(2,2)}=\frac{1}{h_1h_2}\,
\begin{bmatrix}
\phantom{-}1&-2&\phantom{-}1\\-2&\phantom{-}4&-2\\\phantom{-}1&-2&\phantom{-}1
\end{bmatrix}.
\end{equation}
Section~\SEC{numexp} presents numerical results using this
discretization.

%-----------------------[FASTEV]-----------------------------
\section{Fast evaluation of discrete transforms}
\label{sec:fastev}
Consider the evaluation of the grid~$\vctr{h}$ discrete (sub)transform: 
$\forall\vctr{x}\in\{\vctr{x}_{\vctr{i}}\}$,
\begin{equation}\label{eq:subtrans}
S^{\vctr{h},\vctr{l}}(\vctr{x})=\sum_{\vctr{j}=0}^{\vctr{n}} G^{\vctr{l}}(\vctr{x},\vctr{y}_{\vctr{j}}) U_{\vctr{j}}^{\vctr{h},\vctr{l}}.
\end{equation}
The kernel~$G^{\vctr{l}}(\vctr{x},\vctr{y})$ is assumed to be 
asymptotically smooth in~$x_k$ and~$y_k$.
This implies that~$G^{\vctr{l}}(\vctr{x},\vctr{y})$ 
is increasingly smooth in~$x_k$ and~$y_k$,
so that for all allowed errors,~$\epsilon\in\IR$,~$\epsilon>0$, 
and all scales~$H\in\IR$,~$H>0$, there exist~$p,m\in\IN$ for which
a {\em softened kernel}~$G^{\vctr{l}}_{H\vctr{e}_k}(\vctr{x},\vctr{y})$ can be defined 
with the properties:
\begin{itemize}
\item[(\/{\em i}\/)] {\textit{Locality:}}
$G^{\vctr{l}}_{H\vctr{e}_k}(\vctr{x},\vctr{y})=G^{\vctr{l}}(\vctr{x},\vctr{y})$ for 
$|y_k-x_k|\geq mH$.
\item[(\/{\em ii}\/)]
$G^{\vctr{l}}_{H\vctr{e}_k}(\vctr{x},\vctr{y})$ is {\textit{suitably smooth}} in~$x_k$ and~$y_k$ 
on the scale~$H$.
\end{itemize}
Property (\/{\em ii}\/) means
that, both as a function of~$x_k$ for any fixed~$(\vctr{x}-x_k\vctr{e}_k,y)$
and as a function of~$y_k$ for any fixed~$(x,\vctr{y}-y_k\vctr{e}_k)$, 
$G^{\vctr{l}}_{H\vctr{e}_k}(\vctr{x},\vctr{y})$ can be approximated up to an error of at most~$\epsilon$ by a 
$p$-order interpolation from values 
$\{G^{\vctr{l}}_{H\vctr{e}_k}(\vctr{x}+(jh-x_k)\vctr{e}_k,\vctr{y})\mid{}j\in\IZ\}$ and
$\{G^{\vctr{l}}_{H\vctr{e}_k}(\vctr{x},\vctr{y}+(jh-y_k)\vctr{e}_k)\mid{}j\in\IZ\}$, for any~$h\in(0,H]$.
The order of interpolation~$p$ is required to increase only moderately with 
decreasing~$\epsilon$, in particular, 
$p=O(\mathrm{log}\,1/\epsilon)$ as~$\epsilon\to{0}$.
Suitable smoothness of~$G^{\vctr{l}}_{H\vctr{e}_k}(\vctr{x},\vctr{y})$ 
translates into the requirement that
\begin{equation}
(\gamma_2 H)^{p}
\big|G_{H\vctr{e}_k}^{\vctr{l}(p\vctr{e}_k)}(\vctr{x},\vctr{y})\big|
\leq{O}(\epsilon),
\end{equation}
for any~$(\vctr{x},\vctr{y})$ in the domain of interest. 
Here,~$G^{\vctr{l}(p\vctr{e}_k)}_{H\vctr{e}_k}(\vctr{x},\vctr{y})$ 
denotes a~$p$-order derivative 
with respect to either~$x_k$ or~$y_k$. Further,~$\gamma_2$ is a constant 
depending on the interpolation geometry. In particular,~$\gamma_2=1/2$ for 
the usual central interpolations.

The transform~\EQ{subtrans} can now be replaced by a softened transform and a 
correction: 
\begin{equation}\label{eq:S+M}
S^{\vctr{h},\vctr{l}}(\vctr{x})=S_{H\vctr{e}_k}^{\vctr{h},\vctr{l}}(\vctr{x})+M^{\vctr{h},\vctr{l}}_{H\vctr{e}_k}(\vctr{x}),
\end{equation}
with
\begin{equation}\label{eq:softrans}
S^{\vctr{h},\vctr{l}}_{H\vctr{e}_k}(\vctr{x})=\sum_{\vctr{j}=0}^\vctr{n}G^{\vctr{l}}_{H\vctr{e}_k}(\vctr{x},\vctr{y}_{\vctr{j}})U^{\vctr{h},\vctr{l}}_\vctr{j}
\end{equation}
and
\begin{equation}
\label{eq:correct}
M^{\vctr{h},\vctr{l}}_{H\vctr{e}_k}(\vctr{x})=\sum_{\vctr{j}=0}^\vctr{n}
\left(G^\vctr{l}(\vctr{x},\vctr{y}_{\vctr{j}})-G^{\vctr{l}}_{H\vctr{e}_k}(\vctr{x},\vctr{y}_{\vctr{j}})\right)U^{\vctr{h},\vctr{l}}_\vctr{j}.
\end{equation}
Note that by (\/{\em i}\/), 
the correction~\EQ{correct} is local in the~$k$-direction and
only involves points with~$|{y_k}_{\vctr{j}}-x_k|<mH$.

Exploiting the suitable smoothness of the softened kernel, at 
the expense of an error~${O}(\epsilon)$, one may replace 
$G^{\vctr{l}}_{H\vctr{e}_k}(\vctr{x},\vctr{y}_{\vctr{j}})$ by a~$p$-order interpolation from its values on a grid
$\{\vctr{Y}_{\vctr{J}}\}$ with 
mesh sizes~$\vctr{H}=\vctr{h}+(H-h_k)\vctr{e}_k$.
Specifically, there are interpolation 
weights~$w^{\vctr{h}\vctr{H}}_{\vctr{j}\vctr{J}}$ such that for all~$\vctr{x}$:
\begin{equation}\label{eq:intsum}
S_{H\vctr{e}_k}^{\vctr{h},\vctr{l}}(\vctr{x})=\sum_{\vctr{j}=0}^\vctr{n}\sum_{\vctr{J}\in\Gamma^p_\vctr{j}}
w^{hH}_{\vctr{j}\vctr{J}} G_{H\vctr{e}_k}^\vctr{l}(\vctr{x},\vctr{Y}_{\vctr{J}}) 
U_\vctr{j}^{\vctr{h},\vctr{l}}+{O}(\epsilon).
\end{equation}
where~$\Gamma^p_{\vctr{j}}$ stands for a set of~$p$ nodes in the neighborhood of 
$\vctr{y}_{\vctr{j}}$, \mbox{e.g.}, for central~$p$-order interpolations, 
$\Gamma^p_{\vctr{j}}=\{\vctr{J}\in \IZ \mid 
|{y_i}_{\vctr{j}}-{Y_i}_{\vctr{J}}|\leq{\delta_{ik}pH}/2\}$, with~$\delta_{ik}$ the Kronecker delta. 
Changing the order of summation in~\EQ{intsum} and neglecting 
$O(\epsilon)$ errors:
\begin{equation}\label{eq:fgt}
S_{H\vctr{e}_k}^{\vctr{h},\vctr{l}}(\vctr{x})=\sum_{\vctr{J}=0}^{N} G_{H\vctr{e}_k}^\vctr{l}(\vctr{x},\vctr{Y}_{\vctr{J}}) U_\vctr{J}^{\vctr{H},\vctr{l}},
\end{equation}
with
\begin{equation}\label{eq:cgU}
U_{\vctr{J}}^{\vctr{H},\vctr{l}}=\sum_{\vctr{j}\in{\Gamma^p_\vctr{J}}} w^{\vctr{h}\vctr{H}}_{\vctr{j}\vctr{J}} 
U_\vctr{j}^{\vctr{h},\vctr{l}}.
\end{equation}
The operation~\EQ{cgU} is commonly referred to as anterpolation, since it is 
the {\em adjoint} of {\em interpolation}. 

Next, let~$\{\vctr{X}_{\vctr{I}}\}$ denote an evaluation grid with mesh sizes~$\vctr{H}$. By the 
smoothness of~$G_{H\vctr{e}_k}^{\vctr{l}}(\vctr{x},\vctr{y})$ with respect to~$x_k$ on the scale~$H$, for any 
$\vctr{x}_{\vctr{i}}$ there are interpolation weights~$\bar{w}^{\vctr{h}\vctr{H}}_{\vctr{i}\vctr{I}}$ 
such that for all~$\vctr{y}$:
\begin{equation}\label{eq:intx}
G_{H\vctr{e}_k}^\vctr{l}(\vctr{x}_{\vctr{i}},\vctr{y})=\sum_{\vctr{I}\in{\Gamma^p_{\vctr{i}}}}
\bar{w}^{\vctr{h}\vctr{H}}_{\vctr{i}\vctr{I}} G_{H\vctr{e}_k}^{\vctr{l}}(\vctr{X}_{\vctr{I}},\vctr{y})+
O(\epsilon)
\end{equation}
From~\EQ{intx} it follows that, neglecting~$O(\epsilon)$ errors:
\begin{equation}\label{eq:fgti}
S_{H\vctr{e}_k}^{\vctr{h},\vctr{l}}(\vctr{x}_{\vctr{i}})=
\sum_{\vctr{I}\in{\Gamma^p_{\vctr{i}}}} \bar{w}_{\vctr{i}\vctr{I}}^{\vctr{h}\vctr{H}} S_{H\vctr{e}_k}^{\vctr{H},\vctr{l}}(\vctr{X}_{\vctr{I}}),
\end{equation}
with~$S_{H\vctr{e}_k}^{\vctr{H},\vctr{l}}$ denoting a 
grid~$\vctr{H}$ transform defined as:
\begin{equation}\label{eq:cgt}
S_{H\vctr{e}_k}^{\vctr{H},\vctr{l}}(\vctr{X}_{\vctr{I}})=\sum_{\vctr{J}=0}^{\vctr{N}} G_{H\vctr{e}_k}^\vctr{l}(\vctr{X}_{\vctr{I}},\vctr{Y}_{\vctr{J}}) U_{\vctr{J}}^{\vctr{H},\vctr{l}}
\end{equation}

Summarizing, by~\EQ{fgt} to~\EQ{cgt}, at the expense of an error 
$O(\epsilon)$, the grid~$\vctr{h}$ evaluation can be replaced by:
\begin{enumerate}
\item[(\/{\em i}\/)]
\textit{anterpolation} of~$U^{\vctr{h},\vctr{l}}_{\vctr{j}}$ from the integration grid 
$\{\vctr{y}_{\vctr{j}}\}$ to the coarse integration grid,~$\{\vctr{Y}_{\vctr{J}}\}$, by~\EQ{cgU}.
\item[(\/{\em ii}\/)]
\textit{evaluation} of~$S_{H\vctr{e}_k}^{\vctr{H},\vctr{l}}$ on grid~$\vctr{H}$.
\item[(\/{\em iii}\/)]
\textit{interpolation} of~$S_{H\vctr{e}_k}^{\vctr{H},\vctr{l}}$ from the coarse evaluation grid, 
$\{\vctr{X}_{\vctr{I}}\}$, to the evaluation grid~$\{\vctr{x}_{\vctr{i}}\}$.
\item[(\/{\em iv}\/)] addition of the local \textit{correction}
$M^{\vctr{h},\vctr{l}}_{H\vctr{e}_k}$, for all points of~$\{\vctr{x}_{\vctr{i}}\}$.
\end{enumerate}

Denoting by~$n$ the number of point on grid~$\vctr{h}$, the cost of the 
transfer operations (anterpolation and interpolation) is~$O(pn)$. 
The corrections~\EQ{correct} arise in regions where the kernel is 
insufficiently smooth to be accurately 
approximated by a~$p$-order interpolation from a grid with mesh size
$\vctr{H}$. These regions are of dimension~$d-1$ and the work invested
in the corrections is~$O(mn^{2-1/d})$. It is important to notice that
for~$d=1$, the grid~$h$ evaluation can be transfered to grid~$H$
by~$O(n)$ operations. Of course, the coarse grid evaluation (\/{\em ii}\/) 
can again be replaced by successive (\/{\em i}\/)-(\/{\em iv}\/) to
transfer the multisummation to an even coarser grid. Hence, the 
process can be repeated recursively until a grid is reached at which the 
evaluation can be performed in~$O(n)$ operations by direct summation. 
The grid~$h$ multisummation can thus be evaluated with asymptotically
optimal efficiency, i.e. in~$O(n)$ operations.

If~$d\geq{2}$, however, straightforward evaluation of the corrections
inhibits optimal efficiency. To recover optimal efficiency, it is necessary
to reduce the cost of the corrections to~$O(n)$ operations. 
For this purpose, the correction~\EQ{correct} is rewritten as
\begin{multline}
\label{eq:correct2}
M^{\vctr{h},\vctr{l}}_{H\vctr{e}_k}(\vctr{x})
=
\sum_{\vctr{j}=0}^\vctr{n} 
\Big(
G^{\vctr{l}}_{H\vctr{e}_q}(\vctr{x},\vctr{y}_{\vctr{j}})
-
G^{\vctr{l}}_{H(\vctr{e}_k+\vctr{e}_q)}(\vctr{x},\vctr{y}_{\vctr{j}})
\Big)
U^{\vctr{h},\vctr{l}}_{\vctr{j}}
\\
+
\sum_{\vctr{j}=0}^\vctr{n} 
\Big(
\big(
G^{\vctr{l}}(\vctr{x},\vctr{y}_{\vctr{j}})
-G^{\vctr{l}}_{H\vctr{e}_k}(\vctr{x},\vctr{y}_{\vctr{j}})\big)-
\big(
G^{\vctr{l}}_{H\vctr{e}_q}(\vctr{x},\vctr{y}_{\vctr{j}})-
G^{\vctr{l}}_{H(\vctr{e}_k+\vctr{e}_q)}(\vctr{x},\vctr{y}_{\vctr{j}})
\big)
\Big)
U^{\vctr{h},\vctr{l}}_\vctr{j}\:,
\end{multline}
with~$q\neq{k}$. 
Assuming that a softened kernel inherits its asymptotic 
smoothness properties from the original kernel, the softened kernel
$G^{\vctr{l}}_{H(\vctr{e}_k+\vctr{e}_q)}(\vctr{x},\vctr{y}_{\vctr{j}})$
can be constructed such that it is suitably smooth in~$x_k,x_q,y_k$ and~$y_q$ 
on the scale~$H$ and that the second multisummation in~\EQ{correct2} 
is local in the~$k$- and~$q$-directions.
By the suitable smoothness of the softened kernel in~$x_q$ and 
$y_q$ on the scale~$H$, the first multisummation in~\EQ{correct2} can
be transferred to a grid that is coarse in the $q$-direction.
The process of separating a correction into a softened correction,
that can be transferred to a coarser grid, and a
lower-dimensional correction, can be repeated recursively with respect to
all coordinate directions. The grid~$\vctr{h}$ 
corrections can then be evaluated in~$O(n)$ operations.

Separating the transform~\EQ{cgt} as
\begin{multline}
\label{eq:cgt2}
S_{H\vctr{e}_k}^{\vctr{H},\vctr{l}}(\vctr{X}_{\vctr{I}})=
\sum_{\vctr{J}=0}^{\vctr{N}} 
G_{H(\vctr{e}_k+\vctr{e}_q)}^\vctr{l}(\vctr{X}_{\vctr{I}},\vctr{Y}_{\vctr{J}})
U_{\vctr{J}}^{\vctr{H},\vctr{l}}+
\\
\sum_{\vctr{J}=0}^{\vctr{N}} 
\big(
G_{H\vctr{e}_k}^\vctr{l}(\vctr{X}_{\vctr{I}},\vctr{Y}_{\vctr{J}}) 
-
G_{H(\vctr{e}_k+\vctr{e}_q)}^\vctr{l}(\vctr{X}_{\vctr{I}},\vctr{Y}_{\vctr{J}})
\big)
U_{\vctr{J}}^{\vctr{H},\vctr{l}}\:,
\end{multline}
it is evident that the multisummation~\EQ{cgt} can be transferred to a grid
that is coarse in the~$q$-direction at the expense of a correction 
that is local in the~$q$-direction. Of course, this process can 
also be repeated recursively with respect to all coordinate directions.

Summarizing, to evaluate 
the grid~$\vctr{h}$ discrete (sub)transform~\EQ{subtrans} fast,
the operations (\/{\em i}\/)-(\/{\em iv}\/) are recursively 
applied to transfer the multisummation to grids
that are increasingly coarse in each direction, until
a grid is reached at which the multisummation can be performed in
$O(n)$ operations by direct summation. All 
corrections that arise are treated in the same manner.
The treatment of the corrections ensures that
the correctional work is~${O}(1)$ operations 
per grid point (of the grid on which the corrections are required).
On sufficiently fine grids the corrections are negligible compared to
the discretization error which is made anyway, so that corrections 
can be avoided at all, i.e.,~$m=0$ can be used on the finest grids;
see [9] and appendix A.  However, at this point it is noted that this 
only applies to the magnitude of the evaluation error.
If the evaluation with a certain accuracy is not a final goal, 
but a subtask in the numerical solution of the 
integral equation, then a minimum softening distance~$m=O(p)$ 
is needed to ensure that the fast evaluation operator 
has the same stability properties as the unigrid evaluation operator for 
highly oscillatory components.
If one only considers the evaluation of the integral transform, however,
corrections on the finest grids are unnecessary.
For large scale computations the work involved 
in the fast evaluation is then only determined by the costs of the transfers
on the finer grids, and the additional cost of the coarsest grid 
multisummation. The 
work estimates for the evaluation of all discrete subtransforms 
as indicated in~\cite{Brandt:1998fk} can then indeed be obtained.

Although it is most efficient to apply the softening and coarsening per 
direction, it is usually more convenient to first soften the kernel
with respect to all coordinate directions and then transfer the
multisummation. The additional expenses are only marginal.

\pagebreak

%-----------------------[KERNSOF]-----------------------------
\section{Kernel softening}
\label{sec:kernsof}
In the previous section, we showed that the discrete subtransforms
$S^{\vctr{h},\vctr{l}}(\vctr{x})$, 
resulting from the discretization of~\EQ{int}, 
in principle can be evaluated fast by separating each of the transforms 
in a softened transform and a local correction. The multisummation that 
is required to evaluate the softened transform can then be transferred to a 
coarser grid. As a result of the suitable smoothness of the 
{\em softened kernel} on the coarse grid scale, the evaluation error thus 
introduced is less than the fine grid discretization error. 

In~\cite{Brandt:1998fk} it was shown that for 1-dimensional kernels a convenient 
softening can be obtained by locally replacing the original kernel 
with a polynomial, 
$P_H(\vctr{x},\vctr{y})=\sum_{{i}=0}^{{i}=2p-1}a_{{i}}(y-x)^{{i}}$, 
in such a manner that the resulting kernel is~$p-1$ times continuously
differentiable.
This approach can be extended to multidimensional kernels, 
by allowing the polynomial coefficients to depend on a reduced set of 
variables. In particular, for properly chosen softening distance~$m$
and softening order~$p$,
\begin{equation}
\label{eq:soft1}
G^{l}_{H\vctr{e}_k}(\vctr{t})=\left\{
\begin{array}{ll}
P^l_{H\vctr{e}_k}(\vctr{t})\equiv\sum_{{i}=0}^{2p-1}a_{{i}}(\vctr{t}-t_k\vctr{e}_k)\,{t}_k^{{i}},
&
\qquad|{t}_k|\leq{}mH\\
G^l(\vctr{t}),
&
\qquad\mbox{otherwise},
\end{array}
\right.
\end{equation}
with
$\vctr{t}=\vctr{y}-\vctr{x}$,
defines a softened kernel that is suitably smooth in~${t}_k$ on the scale 
$H$, provided that the coefficients,~$a_{i}(\vctr{t}-t_k\vctr{e}_k)$, satisfy
the continuity conditions
\begin{equation}
\label{eq:contcon}
\sum_{{i}={j}}^{2p-1}a_{{i}}(\vctr{t}-t_k\vctr{e}_k)
\frac{{i}{!}}{({i}-{j}){!}}(\pm{mH})^{{i}-{j}}
=
G^{\vctr{l}-{j}\vctr{e}_k}(\vctr{t}-(t_k\pm{m}H)\vctr{e}_k),
\qquad
j=0,\ldots,p-1.
\end{equation}
The~$2p$ coefficients~$a_{{i}}(\vctr{t}-t_k\vctr{e}_k)$ in 
equation~\EQ{soft1} are uniquely determined by the~$2p$ continuity 
conditions~\EQ{contcon}. One may note that by~\EQ{soft1}, the operation 
is local in~$t_k$. Commonly,~$G^{\vctr{l}}(\vctr{t})$ 
is either an even or an odd 
function of~${t}_k$ and~$a_{{i}}(\vctr{t}-t_k\vctr{e}_k)=0$ for all odd 
${i}$ or all even~${i}$, respectively. Moreover, one should anticipate that
usually the softening polynomial can be condensed to a convenient form
that can be evaluated efficiently.

By~\EQ{contcon}, the polynomial coefficients,~$a_{i}(\vctr{t}-t_k\vctr{e}_k)$,
are a linear combination of the kernel derivatives
$G^{\vctr{l}-{j}\vctr{e}_k}(\vctr{t}-(t_k\pm{m}H)\vctr{e}_k)$, 
$j=0,\ldots,p-1$. 
Therefore, if the kernel~$G^{\vctr{l}}(\vctr{t})$ 
consists of a summation of components, then each of these 
components can be softened independently to form the 
softened kernel. Moreover, the asymptotic smoothness 
properties of the original kernel 
are maintained during the softening operation.
Hence, if the original 
kernel is asymptotically smooth in~$t_q$ (\/$q\neq{k}$\/),
then~$G^{\vctr{l}}_{He_k}(\vctr{t})$ can be softened in the \mbox{$q$-direction}
to create a kernel that is suitably smooth in~${t}_k$ and~${t}_{q}$.
The resulting kernel again inherits its asymptotic smoothness 
properties from the original kernel.
Consequently, if the original kernel
is asymptotically smooth in~$\vctr{t}$, then sequential application of the
softening operation with respect to each coordinate direction yields
a softened kernel that is suitably smooth in~$\vctr{t}$. 

To illustrate the multidimensional softening procedure, we consider the
softening of the kernel~\EQ{G22}. Notice that the kernel consists of
a sum of components: 
$G^{(2,2)}(\vctr{t})=\sum_{{i}=0}^{{i}=5}G_{{i}}^{(2,2)}(\vctr{t})$,
with
%------
\begin{alignat}{2}
G^{(2,2)}_0(\vctr{t})&=\frac{1}{2} {t}_{1} |{t}_{1}| {t}_{2}\,
\mbox{arcsinh}\left(\frac{{t}_2}{{t}_1}\right)\:,
&
\qquad
G^{(2,2)}_3({\vctr{t}})&=\frac{1}{2} {t}_{2} |{t}_{2}| {t}_{1}\,
\mbox{arcsinh}\left(\frac{{t}_1}{{t}_2}\right)\:,
\nonumber\\
G^{(2,2)}_1(\vctr{t})&=-{\frac{1}{6}}{t}_{1}^2 \sqrt{{t}_{1}^{2}+{t}_{2}^{2}},
&
\qquad
G^{(2,2)}_4({\vctr{t}})&=-{\frac{1}{6}}{t}_2^2 \sqrt{{t}_1^2+{t}_2^2},
\label{eq:kerncom}\\
G^{(2,2)}_2({\vctr{t}})&={\frac{1}{6}} |{t}_1^3|,
&
\qquad
G^{(2,2)}_5({\vctr{t}})&={\frac{1}{6}}|{t}_2^{3}|.
\nonumber
\end{alignat}
The component~$G^{(2,2)}_5(\vctr{t})$ requires no softening with 
respect to~${t}_1$ since it is already sufficiently smooth. 
Assuming that identical softening parameters (\/$m$ and~$p$\/) 
are chosen in both coordinate directions,
in regions where softening with respect to~${t}_1$ is required, 
i.e. for~$|{t}_1|\leq{}mH$, the softening polynomials read
\begin{align}
P^{(2,2)}_{0\,H\vctr{e}_1}(\vctr{t})&=
\frac{1}{2} {t}_{1}^{2} {t}_{2}\,
\mbox{arcsinh}\left(\frac{{t}_2}{mH}\right)\\
&+\frac{(mH)^{2p}}
{\Bigl( (mH)^{2} + {t}_{2}^{2} \Bigr)^{p-{\frac{3}{2}}}}
\sum_{{i}=1}^{p-1}
\sum_{{j}=1}^{p}A^{(0)}_{{i}{j}}
\left(\frac{{t}_1}{mH}\right)^{2{i}}
\left(\frac{{t}_2}{mH}\right)^{2{j}}
\label{eq:s0}\\
%-------------------%
P^{(2,2)}_{1\,H\vctr{e}_1}(\vctr{t})&=
\frac{(mH)^{2p}}
{\Bigl( (mH)^{2} + {t}_{2}^{2} \Bigr)^{p-{\frac{3}{2}}}}
\sum_{{i}=1}^{p-1}
\sum_{{j}=1}^{p}A^{(1)}_{{i}{j}}
\left(\frac{{t}_1}{mH}\right)^{2{i}}
\left(\frac{{t}_2}{mH}\right)^{2{j}}
\label{eq:s1}\\
%-------------------%
P^{(2,2)}_{2\,H\vctr{e}_1}(\vctr{t})&=
(mH)^{3}
\sum_{{i}=1}^{p-1}
A^{(2)}_{{i}}
\left(\frac{{t}_1}{mH}\right)^{2{i}}
\label{eq:s2}\\
%-------------------%
P^{(2,2)}_{3\,H\vctr{e}_1}(\vctr{t})&=
\frac{(mH)^{2p}}
{\Bigl( (mH)^{2} + {t}_{2}^{2} \Bigr)^{p-{\frac{3}{2}}}}
\sum_{{i}=1}^{p-1}
\sum_{{j}=1}^{p}A^{(3)}_{{i}{j}}
\left(\frac{{t}_1}{mH}\right)^{2{i}}
\left(\frac{{t}_2}{mH}\right)^{2{j}}
\nonumber\\ 
&+\mbox{arcsinh}\left|\frac{mH}{{t}_2}\right|
\frac{(mH)^{2p-1}}
{\Bigl( (mH)^{2} + {t}_{2}^{2} \Bigr)^{p-2}}
\sum_{{i}=1}^{p-1}
\sum_{{j}=1}^{p}B^{(3)}_{{i}{j}}
\left(\frac{{t}_1}{mH}\right)^{2{i}}
\left(\frac{{t}_2}{mH}\right)^{2{j}}
\label{eq:s3}
\end{align}
\begin{align}
%-------------------%
P^{(2,2)}_{4\,H\vctr{e}_1}(\vctr{t})&=
\frac{(mH)^{2p}}
{\Bigl( (mH)^{2} + {t}_{2}^{2} \Bigr)^{p-{\frac{3}{2}}}}
\sum_{{i}=1}^{p-1}
\sum_{{j}=1}^{p}A^{(4)}_{{i}{j}}
\left(\frac{{t}_1}{mH}\right)^{2{i}}
\left(\frac{{t}_2}{mH}\right)^{2{j}}
\label{eq:s4}
\end{align}
%-------------------%
As an example, the coefficients in equations~\EQ{s0} to~\EQ{s4} are listed 
in Table~\TAB{coeff4} for softening-order~$p=4$. Note that by~\EQ{soft1}, 
$G^{(2,2)}_{{i}\,H\vctr{e}_1}(\vctr{t})=G^{(2,2)}_{{i}}(\vctr{t})$ if~$|{t}_1|>mH$.

%A4-------------------------------------------
\begin{table}
\begin{tabular}{|c|c|c|}
\cline{1-3}
&&\\
$
A^{(0)}=
{\left(\begin {array}{ccccc} 
0&{\frac {5}{32}}&{\frac {5}{24}}&{\frac{1}{12}}&0
\\\noalign{\medskip}
0&-{\frac{1}{32}}&{\frac{3}{16}}&{\frac{1}{8}}&0
\\\noalign{\medskip}
0&-{\frac {5}{32}}&-{\frac{1}{2}}&-{\frac{1}{4}}&0
\\\noalign{\medskip}
0&{\frac{1}{32}}&{\frac {5}{48}}&{\frac{1}{24}}&0
\end {array}\right)}
$
&
%B4-------------------------------------------
$
A^{(1)}=
{\left (\begin {array}{ccccc} 
{\frac {1}{96}}&0&0&0&0
\\\noalign{\medskip}
-{\frac {3}{32}}&-{\frac {5}{16}}&-{\frac {5}{12}}&-{\frac{1}{6}}&0
\\\noalign{\medskip}
-{\frac {3}{32}}&-{\frac {5}{24}}&-{\frac{1}{12}}&0&0
\\\noalign{\medskip}
{\frac {1}{96}}&{\frac{1}{48}}&0&0&0
\end {array}\right )}
$
&
%C14--------------------------------------------
$
A^{(2)}=
{\left( \begin{array}{cccc}
       -{\frac {1}{96}}&{\frac {3}{32}}&
        {\frac {3}{32}}&-{\frac {1}{96}}
\end{array} \right)}
$
\\
&&\\
\cline{1-3}
&&\\
%C24--------------------------------------------
$
A^{(3)}=
{\left (\begin {array}{ccccc} 
0&-{\frac {23}{96}}&-{\frac {35}{96}}&-{\frac {5}{32}}&0
\\\noalign{\medskip}
0&{\frac {7}{32}}&{\frac {5}{32}}&{\frac{1}{32}}&0
\\\noalign{\medskip}
0&{\frac{1}{32}}&{\frac {9}{32}}&{\frac {5}{32}}&0
\\\noalign{\medskip}
0&-{\frac {1}{96}}&-{\frac {7}{96}}&-{\frac{1}{32}}&0
\end {array}\right )}
$
&
%C24--------------------------------------------
$
B^{(3)}=
{\left (\begin {array}{ccccc} 
0&{\frac {5}{32}}&{\frac {5}{16}}&{\frac {5}{32}}&0
\\\noalign{\medskip}
0&{\frac {15}{32}}&{\frac {15}{16}}&{\frac {15}{32}}&0
\\\noalign{\medskip}
0&-{\frac {5}{32}}&-{\frac {5}{16}}&-{\frac {5}{32}}&0
\\\noalign{\medskip}
0&{\frac{1}{32}}&{\frac{1}{16}}&{\frac{1}{32}}&0
\end {array}\right )}
$
&
%D4------------------------------------
$
A^{(4)}=
{\left (\begin {array}{ccccc} 
0&-{\frac {5}{96}}&-{\frac {5}{16}}&-{\frac {5}{12}}&-{\frac{1}{6}}
\\\noalign{\medskip}
0&-{\frac {5}{32}}&-{\frac {5}{24}}&-{\frac{1}{12}}&0
\\\noalign{\medskip}
0&{\frac {5}{96}}&{\frac{1}{48}}&0&0
\\\noalign{\medskip}
0&-{\frac {1}{96}}&0&0&0\end {array}\right )}
$
\\
&&\\
\cline{1-3}
%E4-------------------------------------------
\multicolumn{3}{|c|}{}\\
\multicolumn{3}{|c|}{$
A^{(5)}={{
{\left (\begin {array}{cccc} {\frac {161\,\sqrt {2}}{6144}}&{\frac {\sqrt {2}\left (-
671+480\,\sqrt {2}\ln (1+\sqrt {2})\right )}{6144}}&-{\frac {517\,\sqrt {2}}{6144}}&{
\frac {67\,\sqrt {2}}{6144}}\\\noalign{\medskip}{\frac {\sqrt {2}\left (-671+480\,
\sqrt {2}\ln (1+\sqrt {2})\right )}{6144}}&{\frac {\sqrt {2}\left (321+2880\,\sqrt {2}
\ln (1+\sqrt {2})\right )}{6144}}&{\frac {\sqrt {2}\left (-480\,\sqrt {2}\ln (1+\sqrt {
2})-261\right )}{6144}}&{\frac {\sqrt {2}\left (96\,\sqrt {2}+35\right )}{6144}}
\\\noalign{\medskip}-{\frac {517\,\sqrt {2}}{6144}}&{\frac {\sqrt {2}\left (-480\,
\sqrt {2}\ln (1+\sqrt {2})-261\right )}{6144}}&{\frac {83\,\sqrt {2}}{2048}}&-{\frac {
47\,\sqrt {2}}{6144}}\\\noalign{\medskip}{\frac {67\,\sqrt {2}}{6144}}&{\frac {\sqrt 
{2}\left (96\,\sqrt {2}+35\right )}{6144}}&-{\frac {47\,\sqrt {2}}{6144}}&{\frac {3\,
\sqrt {2}}{2048}}\end {array}\right )}}}
$}\\
\multicolumn{3}{|c|}{}\\
\cline{1-3}
\end{tabular}
\caption{Coefficients in equations~\EQ{s0} to~\EQ{doublesoft}
for softening-order~$p=4$.
\label{tab:coeff4}}
\end{table}

To obtain a softened kernel that is smooth in both~${t}_1$ and~${t}_2$, 
subsequently, the kernel is softened with respect to~${t}_2$.
Because the original kernel has the symmetry property  
$G^{(2,2)}({t}_1,{t}_2)=G^{(2,2)}({t}_2,{t}_1)$, equations~\EQ{s0} to~\EQ{s4} 
determine the softened kernel in regions where the original kernel is 
suitably smooth in either~${t}_1$ or~${t}_2$. To obtain the softened
kernel in the region where the original kernel is unsmooth in both 
${t}_1$ and~${t}_2$, the softening operation with respect to~${t}_2$ 
is applied to the softened kernel~$G^{(2,2)}_{H\vctr{e}_1}(\vctr{t})$. 
For~$\vctr{t}\in[-mH,mH]^2$, this yields the polynomial
\begin{multline}
\label{eq:doublesoft}
P^{(2,2)}_{(H,H)}(\vctr{t})=
(mH)^3
\sum_{{i}=1}^{p-1}
\sum_{{j}=1}^{p-1}A^{(5)}_{{i}{j}}
\Big(\frac{{t}_1}{mH}\Big)^{2{i}}
\Big(\frac{{t}_2}{mH}\Big)^{2{j}}
+\\
\sum_{{i}=1}^{p-1}
A^{(2)}_{{i}}
\bigg(\Big(\frac{{t}_1}{mH}\Big)^{2{i}}+
\Big(\frac{{t}_2}{mH}\Big)^{2{i}}\bigg).
\end{multline}
For~$p=4$ the coefficients in~\EQ{doublesoft} are listed in 
Table~\TAB{coeff4}. One may note that~$A^{(5)}_{{i}{j}}$ is 
symmetric, so that the symmetry  of the original kernel in~${t}_1$ 
and~${t}_2$ is maintained. 

\begin{figure}
\begin{center}
%\psfrag{t_1}{\small{\raisebox{-0pt}{\hspace{0pt}\small{$t_1$}}}}
%\psfrag{t_2}{\small{\raisebox{-0pt}{\hspace{0pt}\small{$t_2$}}}}
%\psfrag{mH}{\small{\raisebox{-0pt}{\hspace{0pt}\small{$mH$}}}}
%\psfrag{Theta}{\small{\raisebox{-0pt}{\hspace{0pt}\small{}}}}
%\psfrag{1}{\small{\raisebox{-0pt}{\hspace{0pt}\small{\ding{192}}}}}
%\psfrag{2}{\small{\raisebox{-0pt}{\hspace{0pt}\small{\ding{193}}}}}
%\psfrag{3}{\small{\raisebox{-0pt}{\hspace{0pt}\small{\ding{194}}}}}
%\psfrag{4}{\small{\raisebox{-0pt}{\hspace{0pt}\small{\ding{195}}}}}
\includegraphics[width=0.5\textwidth]{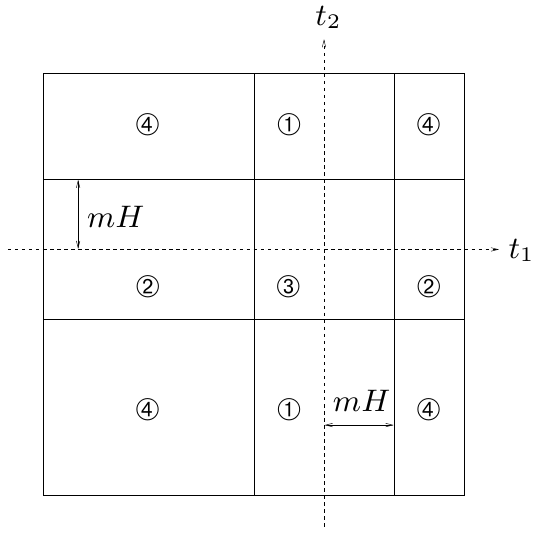}
\caption{Division of domain in softening regions: 
\ding{192}: softening in~$t_1$ only, 
\ding{193}: softening in~$t_2$ only, 
\ding{194}: softening in~$t_1$ and~$t_2$,
\ding{195}: original kernel is sufficiently smooth.
\label{fig:fig1}}
\end{center}
\end{figure}

From equations~\EQ{s0} to~\EQ{doublesoft} 
it follows that the softened kernel is given by
\begin{equation}
\label{eq:g22soft}
G^{(2,2)}_{(H,H)}(\vctr{t})=
\left\{ 
\begin{array}{lll}
G^{(2,2)}_{5}({t}_1,{t}_2)+
\displaystyle{\sum_{{i}=0}^{4}}
P^{(2,2)}_{{i}\,H\vctr{e}_1}({t}_1,{t}_2)
\quad&
|{t}_1|\leq{mH},|{t}_2|>mH
\qquad&
\mbox{\ding{192}}
\\
G^{(2,2)}_{5}({t}_2,{t}_1)+
\displaystyle{\sum_{{i}=0}^{4}}
P^{(2,2)}_{{i}\,H\vctr{e}_1}({t}_2,{t}_1)&
|{t}_1|>mH,|{t}_2|\leq{}mH&
\mbox{\ding{193}}
\\
P^{(2,2)}_{(H,H)}(\vctr{t})&
|{t}_1|\leq{}mH,|{t}_2|\leq{}mH&
\mbox{\ding{194}}
\\[3mm]
G^{2,2}(\vctr{t})&
\mbox{otherwise}.&
\mbox{\ding{195}}
\end{array}
\right.
\end{equation}
The encircled numbers in~\EQ{g22soft} refer to Figure~\FIG{fig1}.

%-----------------------[WOMIN]-----------------------------
\section{Work minimization and error control}
\label{sec:womin}
The multilevel fast evaluation algorithm relies on the smoothness of the 
softened kernel. For the softened kernels resulting from the operation 
defined in Section~\SEC{kernsof}, 
the smoothness depends on the softening order,~$p$, 
and the softening distance,~$m$. The work invested in the fast 
evaluation algorithm also depends on~$p$ and~$m$.
 In this section we derive the~$m$ and~$p$ for which the computational work is 
minimized subject to the condition that the incremental evaluation error does
not exceed the discretization error. 

To obtain the constraint for~$\vctr{m}$ and~$\vctr{p}$, 
we analyze the additional 
evaluation error on the target grid,~$\vctr{h}$, 
due to transferring the evaluation 
of the softened transform from a grid with mesh sizes~$\vctr{H}/2$ 
to a grid with mesh sizes~$\vctr{H}$. 
This error results from transferring the evaluation from the 
integration grid~$\{\vctr{y}^{\vctr{H}/2}_{\vctr{j}}\}$ to the grid 
$\{\vctr{y}^\vctr{H}_{\vctr{j}}\}$, 
i.e. from replacing the  softened kernel 
$G_\vctr{H}^{\vctr{l}}(\vctr{x},\vctr{y})$ for fixed 
$\vctr{x}$ and as a function of~$\vctr{y}$ by a~\mbox{$\vctr{p}$-order} 
interpolation from 
$\{\vctr{y}^\vctr{H}_{\vctr{j}}\}$, and from replacing the transforms on 
the evaluation grid~$\{\vctr{x}^{\vctr{H}/2}_{\vctr{i}}\}$ 
by an interpolation of transforms on 
$\{\vctr{x}^\vctr{H}_{\vctr{i}}\}$, {i.e.} 
from replacing~$G_{\vctr{H}}^{\vctr{l}}(\vctr{x},\vctr{y})$ 
for fixed~$\vctr{y}$ and as a function of~$\vctr{x}$ by a~$\vctr{p}$-order 
interpolation from~$\{\vctr{x}^\vctr{H}_\vctr{i}\}$.

The additional evaluation error,~$E(\vctr{x}_{\vctr{i}}^{\vctr{h}})$, is the 
sum of the local 
interpolation errors, that is, the difference between the actual kernel value 
and its interpolation from the grids~$\{\vctr{y}^\vctr{H}_{\vctr{j}}\}$ and 
$\{\vctr{x}^\vctr{H}_{\vctr{i}}\}$, weighted by 
$U^{\vctr{h},\vctr{l}}_{\vctr{j}}$:
\begin{equation}
\label{eq:Ey}
E(\vctr{x}_\vctr{i}^\vctr{h})=
\sum_{\vctr{j}=0}^\vctr{n}
\bigg(
G^{\vctr{l}}_\vctr{H}(\vctr{x}_\vctr{i}^\vctr{h},\vctr{y}^\vctr{h}_\vctr{j})-
\sum_{\vctr{J}\in\Gamma^{\vctr{p}}_\vctr{j}} 
w^{\vctr{h}\vctr{H}}_{\vctr{j}\vctr{J}}
\sum_{\vctr{I}\in\Gamma^\vctr{p}_\vctr{i}} 
\bar{w}^{\vctr{h}\vctr{H}}_{\vctr{i}\vctr{I}} 
G^{\vctr{l}}_\vctr{H}(\vctr{x}^\vctr{H}_\vctr{I},\vctr{y}^\vctr{H}_\vctr{J})\bigg)
U^{\vctr{h},\vctr{l}}_\vctr{j},
\end{equation}
In regions where~$u(\vctr{y})$ is~$\vctr{s}$-times differentiable,
\begin{equation}
\label{eq:UL}
\left|U^{\vctr{h},\vctr{l}}_\vctr{j}\right|\leq
\big\|u^{(\vctr{s})}\big\|_{\text{max},\Omega}
\prod_{k=1}^{d} \left(2\left(\gamma_3 h_k\right)^{s_k-l_k+1}+
O\left(h_k^{s_k-l_k+2}\right)\right),
\end{equation}
with~$\gamma_3\approx{}0.5$ for a uniform grid;
see~\cite{Brandt:1998fk}. Denoting by~$\Theta=\{\vctr{y}-\vctr{x}\mid
\vctr{x}\in\underline\Omega,\vctr{y}\in\Omega\}$, it follows from~\EQ{Ey} 
and~\EQ{UL} that the error per unit of integration
caused by the transfer from grid~$\vctr{H}/2$ to grid~$\vctr{H}$ is bounded by
\begin{equation}
\label{eq:e_L}
|\varepsilon^\vctr{H}|\leq
\alpha_2\,\big\|G^{\vctr{l}}_\vctr{H}- 
\II^{\vctr{H}/2}_{\vctr{H}}G^{\vctr{\vctr{l}}}_\vctr{H}\big\|_{1,\Theta}
\,\big\|u^{({\vctr{s}})}\big\|_{\text{max},\Omega}
\,\prod_{k=1}^{d} \gamma_3 \left( \gamma_3 h_k \right)^{s_k-l_k},
\end{equation}
where~$\alpha_2$ is some positive constant and
$\big\|G^{\vctr{l}}_\vctr{H}- 
\II^{\vctr{H}/2}_{\vctr{H}}G^{\vctr{\vctr{l}}}_\vctr{H}\big\|_{1,\Theta}$
denotes the average 
in~$\Theta$ of the absolute value of the~$\vctr{p}$-order interpolation error,
introduced by replacing~$G^{\vctr{l}}_\vctr{H}(\vctr{x},\vctr{y})$ for 
fixed~$\vctr{x}$ by an 
interpolation from~$\{\vctr{y}^\vctr{H}_\vctr{j}\}$ and for fixed 
$\vctr{y}$ by an interpolation 
from~$\{\vctr{x}^\vctr{H}_\vctr{i}\}$. 
From~\EQ{tau^h} and~\EQ{e_L} it follows that the incremental evaluation error 
is smaller than the fine grid discretization error if the following 
requirement is satisfied: 
\begin{equation}
\label{eq:creq}
\big\|G^{\vctr{l}}_\vctr{H}- 
\II^{\vctr{H}/2}_{\vctr{H}}G^{\vctr{\vctr{l}}}_\vctr{H}\big\|_{1,\Theta}
\prod_{k=1}^{d}\gamma_3\left(\gamma_3 h_k \right)^{s_k-l_k}
\leq
(\alpha_1/\alpha_2)\,
\sum_{k=1}^{d}
(\gamma_1 h_k)^{s_k}
\big\|G\big\|_{1,\Theta}
\frac
{\big\|u^{(s_k{e}_k)}\big\|_{\text{max},\Omega}}
{\big\|u^{(\vctr{s})}\big\|_{\text{max},\Omega}}.
\end{equation}
Notice that in the one dimensional case the derivatives of~$u(\vctr{y})$ in 
the right-hand side of~\EQ{creq} cancel. The relation between the evaluation 
error and the discretization error is then independent of~$u(\vctr{y})$.

The interpolation error is composed of the interpolation error per direction.
Clearly, the requirement that the incremental evaluation error
is smaller than the fine grid discretization error is satisfied if for 
every direction the contribution to the evaluation error is smaller than 
the contribution to the discretization error. Hence, requirement~\EQ{creq} 
can be separated in the following requirement per direction:
\begin{equation}
\label{eq:sreq}
\big\|G^{\vctr{l}}_\vctr{H}- 
\II^{\vctr{H}-H_k\vctr{e}_k/2}_{\vctr{H}}G^{\vctr{\vctr{l}}}_\vctr{H}\big\|_{1,\Theta}
\leq
h_k^{s_k-d(s_k-l_k)}
\,c^{}_{G,u}\,c_\vctr{h}
\,c_{\gamma}
\:,
\end{equation}
where
\begin{equation}
c_{G,u}=\big\|G\big\|_{1,\Theta}
\displaystyle{
\frac
{\big\|u^{(s_k\vctr{e}_k)}\big\|_{\text{max},\Omega}}
{\big\|u^{(\vctr{s})}\big\|_{\text{max},\Omega}}
},
\:\:
c_\vctr{h}=\prod_{i=1}^{d}
\frac
{h_k^{s_k-l_k}}{h_i^{s_i-l_i}},
\:\:
c_{\gamma}=(\alpha_1/\alpha_2)\,\gamma_1^{s_k}
\prod_{i=1}^{d}\gamma_3^{-(s_i-l_i+1)}.
\nonumber
\end{equation}
The average interpolation error on the left-hand side of~\EQ{sreq} 
depends on the properties of~$G^{\vctr{l}}_\vctr{H}(\vctr{x},\vctr{y})$, and 
thus on the choice of~$\vctr{m}$ and~$\vctr{p}$.
The specific dependence of~\EQ{sreq} on~$\vctr{m}$ and 
$\vctr{p}$ is derived in Appendix~A.

The computational work per grid~$\vctr{H}/2$ node involved in transferring 
the evaluation of the discrete transform from grid~$\vctr{H}/2$ to 
grid~$\vctr{H}$ is estimated
\begin{equation}
\label{eq:W_L}
W^\vctr{H}=O\left(2\left[1-2^{-d}\right]p_k+4d\,m_k\bar{H}{}^{1-d}\right),
\end{equation}
where~$\bar{H}^{1-d}$ estimates the cost of the~$d-1$ dimensional summation 
associated with the correction. The estimate~\EQ{W_L}
assumes that~$p_k$ and~$m_k$ are independent of~$k$. This work estimate is 
obtained as follows: defining an operation to be one multiplication and 
one addition, the number of operations involved in the~$\vctr{p}$-order 
anterpolation from a grid with mesh sizes~$\vctr{H}/2$ to a semi-coarse grid 
with mesh sizes~$(\vctr{H}+\vctr{e}_kH_k)/2$ is~$p_k/2$, 
since for half of the values the transfer is trivial. The number of nodes 
on this semi-coarse grid is approximately half the number of nodes on the
$\vctr{H}/2$ grid. Hence, the next step in the 
anterpolation is performed in 
$p_{k}/4$ operations. In general, the~$i^{th}$ step 
in the anterpolation takes~$p_k/2^{i}$ operations. The number of operations 
resulting from the interpolation is obtained in a similar manner, so that 
the total amount of work invested in the transfer operations is 
approximately~$2[1-2^{-d}]p_{k}$. The evaluation of the 
corrections~\EQ{correct} involves~$4m_{k}$ summations over a~$d-1$ 
dimensional domain per direction per grid~$\vctr{H}/2$ node. 
These summations are again evaluated fast, so that the cost per grid point
is~$O(1)$. One should note that
an accurate estimate of the cost of the corrections is not essential,
since~$m_k=0$ is employed on the finest grids; see~\cite{Brandt:1998fk} and Appendix~A.
Assuming that the dimensions of the domains~$\underline{\Omega}$ and~$\Omega$ 
are~$O(1)$,  we arrive at the total work estimate~\EQ{W_L}. 

The optimal transfer parameters are the~$\vctr{m}$ and~$\vctr{p}$ that 
minimize~\EQ{W_L} subject to~\EQ{sreq}. As an example, 
the optimization process 
for the fast evaluation of the model problem is presented in Appendix~B.

%-----------------------[NUMEXP]-----------------------------
\section{Numerical experiments}
\label{sec:numexp}
The algorithm is tested for the 
integral transform with kernel~\EQ{modker} on a domain $\Omega=[-1,1]^2$ with
\begin{equation}\label{e:utest}
u(\vctr{y})=
\left\{
\begin{array}{ll}
\displaystyle{
\prod_{k=1}^{2}
\left[
-\frac{1}{3}+\left(\frac{10\,y_k}{9}\right)^2
-\frac{2}{3}\left|\frac{10\,y_k}{9}\right|^3\right]
}\:,
&
\qquad\vctr{y}\in[-9/10,9/10]^2\:,\\
\\
0,
&
\qquad\mbox{otherwise}.
\end{array}
\right.
\end{equation}
The above problem is solved 
numerically using an~$\vctr{s}=2$ discretization on uniform
grids with mesh widths~$h_k=2^{1-K}$,~$K=5,6,\ldots,11$.
All boundary terms~\EQ{B} vanish 
and only the discrete transform~$S^{{h},(2,2)}(\vctr{x})$ by~\EQ{S22} 
requires evaluation.

To evaluate this transform fast, the softening of 
$G^{(2,2)}(\vctr{x},\vctr{y})$ presented in Section~\SEC{kernsof} is used. 
Details of the derivation of the optimal transfer parameters are presented
in Appendix B.
The parameters used in the 
computations are obtained as follows: 
first~${p}_k^{*}$ is calculated using equation~\EQ{pk}, with~$l_k=2$, 
$h_k=2^{1-K}$,~$H_k=2^{1-L}$ and the constant~$c_a$ 
in equation~\EQ{lng} set to~$0$. Next,~$p_k$ is obtained from:
\begin{equation}
p_k=\left\{
\begin{array}{ll}
2\,\big\lfloor{p}_k^{*}/2+1\big\rfloor,
&
\qquad{p}_k^{*}\geq l_k+2,
\\
l_k+2,
&
\qquad\mbox{otherwise},
\end{array}
\right.
\end{equation}
Finally,~$m_k$ is obtained from:
\begin{equation}
m_k=\left\{
\begin{array}{ll}
\Big\lfloor\half+
\frac{3}{32}\bar{H}\gamma_2^{-1}\chi\left(p_k-l_k-1\right)\big\rfloor,
&
\qquad{p}_k^{*}\geq l_k+2,
\\
0,
&
\qquad\mbox{otherwise},
\end{array}
\right.
\end{equation}
with~$\bar{H}$ set to~$H_k/2$. The values of~$p_k$ and~$m_k$ thus obtained
are listed in Table~\TAB{mp1}. 
The table confirms that~$m_k=0$ can indeed be used on several of the finer 
grids and that for larger~$K$ the number of such grids increases.
\begin{table}
\begin{small}
\caption{Transferparameters~$p_k$ and~$m_k$  used in the evaluation 
of~$S^{\vctr{h},(2,2)}$.
\label{tab:mp1}}
\begin{center}
\begin{tabular}[t]{|c||c|c||c|c||c|c||c|c||c|c||c|c||c|c|}
\cline{2-15}
\multicolumn{1}{c|}{}&
\multicolumn{2}{c||}{$K=5$}&
\multicolumn{2}{c||}{$K=6$}&
\multicolumn{2}{c||}{$K=7$}&
\multicolumn{2}{c||}{$K=8$}&
\multicolumn{2}{c||}{$K=9$}&
\multicolumn{2}{c||}{$K=10$}&
\multicolumn{2}{c|}{$K=11$}\\
\cline{1-15}
$L$&$p_k$&$m_k$&$p_k$&$m_k$&$p_k$&$m_k$&$p_k$&$m_k$&$p_k$&$m_k$&$p_k$&$m_k$&$p_k$&$m_k$\\
\cline{1-15}
K-1&4&0&4&0&4&0&4&0&4&0&4&0&4&0 \\
\cline{1-15}
K-2&6&2&4&1&4&0&4&0&4&0&4&0&4&0 \\
\cline{1-15}
K-3&&&6&3&6&3&6&2&4&1&4&1&4&0 \\
\cline{1-15}
K-4&&&&&8&5&8&4&6&3&6&3&6&2 \\
\cline{1-15}
K-5&&&&&&&10&6&8&5&8&5&8&4 \\
\cline{1-15}
K-6&&&&&&&&&10&8&10&7&10&6 \\
\cline{1-15}
\end{tabular}
\end{center}
\end{small}
\end{table}
%------------------------------------------------------------------

To monitor the accuracy of the fast evaluation in relation to
the discretization 
error, the error~$\epsilon_K^L$ is defined by the~$l_2$-norm of
the difference between the exact solution and 
the numerical solution that is 
obtained on level~$K$ when direct multisummation 
is performed on level~$L$:
\begin{equation}
\epsilon_K^L=
\Bigg(
\frac{1}{n}
\sum_{\vctr{i}=0}^{\vctr{n}}
|(G^h\tilde{u}^{h})^{K,L}(\vctr{x}^\vctr{h}_\vctr{i})-Gu(\vctr{x}^\vctr{h}_\vctr{i})|^2
\bigg)^{1/2}.
\end{equation}
One may note that~$\epsilon_K^K$ is the~$l_2$-norm of the discretization error 
on the level~$K$ grid. Table~\TAB{mlcerr1} lists 
the errors obtained for the model-problem. 
As a side note, we mention that the corrections are evaluated by means of
the multilevel matrix multiplication algorithm; see~\cite{Brandt:1990kl}. 
The leftmost column confirms~$O(h^2)$-convergence of the discretization 
error. The entries marked 
by asterisks denote the results for~$L=K/2$. In this case, the grid on which 
direct multisummation is performed consists of~$O(\sqrt{n})$ nodes. 
The table clearly shows that with the presented fast evaluation algorithm,
the multisummation can be performed on a grid with 
$O(\sqrt{n})$ points at negligible loss of accuracy. 
\begin{table}
\caption{Error~${\epsilon_{K}^{L}}$ for the model problem.
\label{tab:mlcerr1}}
\begin{center}
\begin{tabular}[t]{|r||r|r|r|r|r|r|r|}
\cline{1-8}
K&
\multicolumn{1}{c|}{$L=K$}&
\multicolumn{1}{c|}{$K-1$}&
\multicolumn{1}{c|}{$K-2$}&
\multicolumn{1}{c|}{$K-3$}&
\multicolumn{1}{c|}{$K-4$}&
\multicolumn{1}{c|}{$K-5$}&
\multicolumn{1}{c|}{$K-6$}\\
\cline{1-8}
\cline{1-8}
5&2.01 10$^{-4}$&2.05 10$^{-4}$&3.60 10$^{-4}$&&&&\\
\cline{1-8}
6&5.17 10$^{-5}$&4.95 10$^{-5}$&1.02 10$^{-4}$&$^{\ast}$1.98 10$^{-4}$ &&&\\
\cline{1-8}
7&1.31 10$^{-5}$&1.27 10$^{-5}$&1.08 10$^{-5}$&1.84 10$^{-5}$&3.18 10$^{-5}$&&\\
\cline{1-8}
8&$\approx$3.3 10$^{-6}$&3.24 10$^{-6}$&2.81 10$^{-6}$&3.12 10$^{-6}$&$^{\ast}$4.18 10$^{-6}$&6.64 10$^{-6}$&\\
\cline{1-8}
9&$\approx$8.1 10$^{-7}$&&7.54 10$^{-7}$&1.06 10$^{-6}$&1.38 10$^{-6}$&2.13 10$^{-6}$&3.55 10$^{-6}$\\
\cline{1-8}
10&$\approx$2.0 10$^{-7}$&&&2.18 10$^{-7}$&2.33 10$^{-7}$&$^{\ast}$2.47 10$^{-7}$&3.21 10$^{-7}$\\ 
\cline{1-8}
11&$\approx$5.0 10$^{-8}$&&&&4.23 10$^{-8}$&4.73 10$^{-8}$&5.28 10$^{-8}$\\
\cline{1-8}
\end{tabular}
\end{center}
\end{table}

To get a better insight into the error introduced by the fast evaluation,
we also monitor the {\em incremental} error,
defined by the~$l_2$-norm of the difference in the solution on level
$K$ when direct summation is performed on level~$L+1$ and when direct
summation is performed on level~$L$:
\begin{equation}
I\epsilon^L_K=
\Bigg(
\frac{1}{n}
\sum_{\vctr{i}=0}^{\vctr{n}}
|(G^h\tilde{u}^{h})^{K,L+1}(\vctr{x}^\vctr{h}_\vctr{i})-
(G^h\tilde{u}^{h})^{K,L}(\vctr{x}^\vctr{h}_\vctr{i})|^2
\bigg)^{1/2}.
\end{equation}
This quantity measures the additional error introduced by 
transferring the evaluation from level~$L+1$ to level~$L$. The incremental 
errors are listed in Table~\TAB{mlcinc1}. 
The table shows that the incremental evaluation 
errors are in all relevant cases  of the same order of magnitude as the 
discretization error. 
\begin{table}
\caption{Incremental error~$I\epsilon_K^L$ for the model problem.
\label{tab:mlcinc1}}
\begin{center}
\begin{tabular}[t]{|r||r|r|r|r|r|r|}
\cline{1-7}
K&
\multicolumn{1}{c|}{$L=K-1$}&
\multicolumn{1}{c|}{$K-2$}&
\multicolumn{1}{c|}{$K-3$}&
\multicolumn{1}{c|}{$K-4$}&
\multicolumn{1}{c|}{$K-5$}&
\multicolumn{1}{c|}{$K-6$}\\
\cline{1-7}
\cline{1-7}
5&3.70 10$^{-5}$&1.77 10$^{-4}$&&&&\\
\cline{1-7}
6&3.79 10$^{-6}$&5.94 10$^{-5}$&$^{\ast}$1.06 10$^{-4}$&&&\\
\cline{1-7}
7&5.65 10$^{-7}$&3.55 10$^{-6}$&9.55 10$^{-6}$&1.61 10$^{-5}$&&\\
\cline{1-7}
8&              &5.55 10$^{-7}$&1.10 10$^{-7}$&$^{\ast}$1.95 10$^{-6}$&3.62 10$^{-6}$&\\
\cline{1-7}
9&              &              &4.62 10$^{-7}$&6.17 10$^{-7}$&1.05 10$^{-6}$&2.02 10$^{-6}$\\
\cline{1-7}
10&             &              &              &6.08 10$^{-8}$&$^{\ast}$8.31 10$^{-8}$& 1.22 10$^{-7}$\\
\cline{1-7}
11&             &              &              &              &1.87 10$^{-8}$&2.37 10$^{-8}$\\
\cline{1-7}
\end{tabular}
\end{center}
\end{table}

To determine the computational complexity of the fast evaluation,
the expended operations are counted. The operation-count
is obtained as the sum of the transfer-costs and the
cost of the coarsest grid multisummation (of both the original transform
and the corrections) for all grids involved in the fast evaluation.
For the results in Table~\TAB{mlcerr1}, the computational work per 
grid~$\vctr{h}$ node is displayed in Table~\TAB{mlcw1}. The leftmost 
column of the table shows the costs of direct summation. It can be seen that 
these costs amount to~$O(n^2)$ operations. The entries marked by asterisks are 
the computational costs in case direct multisummation is performed on a grid 
with~$O(\sqrt{n})$ nodes. As expected, the 
costs per grid point decrease for increasing~$K$. It is anticipated
that only the costs of the transfers on the fine grids and of the coarsest 
grid multisummation remain as~$\vctr{h}\to{0}$. The total number 
of operations is then~$W\approx{}2\,\bar{p}+1$. In the present case 
this yields 
$W\approx{}9$. The results in table~\TAB{mlcw1} suggest that this may indeed 
be obtained for sufficiently large~$K$.
\begin{table}
\caption{Work per gridpoint spent in the evaluation of
$S^{\vctr{h},(2,2)}$ for the model problem.
\label{tab:mlcw1}}
\begin{center}
\begin{tabular}{|r||r|r|r|r|r|r|r|r|}
\cline{1-8}
K&L=K&K-1&K-2&K-3&K-4&K-5&K-6\\
\cline{1-8}
\cline{1-8}
5&1089&185&151&&&&\\
\cline{1-8}
6&4225&450&70&$^{\ast}$74&&&\\
\cline{1-8}
7&1.7 10$^{4}$&1368&120&56&71&&\\
\cline{1-8}
8&6.6 10$^{4}$&4743&351&63&$^{\ast}$44&53&\\
\cline{1-8}
9&2.6~$10^{5}$&&1197&97&27&24&27\\
\cline{1-8}
10&1.1~$10^{6}$&&&309&43&$^{\ast}$24&24\\ 
\cline{1-8}
11&4.4~$10^{6}$&&&&94&21&16\\ 
\cline{1-8}
\end{tabular}
\end{center}
\end{table}

%-----------------------[CONCL]-----------------------------
\section{Conclusion}
Motivated by the demand for local grid refinement techniques in practical
applications, this work examined the extension to multiple
dimensions of 
a new algorithm for the fast evaluation of integral transforms 
with asymptotically smooth kernels. The discretization procedure
was outlined. Details were presented for the fast evaluation method
in the instance of multiple dimensions. It was shown that the asymptotic
work estimates in~\cite{Brandt:1998fk} can indeed be obtained, provided that
multilevel evaluation of the corrections is applied. 
The softened kernels in the fast evaluation algorithm were constructed 
by applying the softening operation sequentially
with respect to each coordinate direction. 
The optimization of softening parameters 
for multidimensional transforms was discussed.

The fast evaluation algorithm was tested for a two dimensional model problem.
The results showed that with the new algorithm the evaluation of 
multidimensional transforms is also more efficient than with previous 
algorithms. Moreover, the results confirmed the expected asymptotic
work estimates for the considered test case.

%-----------------------[ACKNOW]-----------------------------
%\section*{Acknowledgment}
%The authors gratefully acknowledge the many helpful suggestions of 
%prof. dr. Achi Brandt of the Weizmann Institute, Rehovot, Israel.

\clearpage
%-----------------------[REFS]-----------------------------
\bibliographystyle{plain}
\bibliography{BibFile}

\clearpage
%-----------------------[APPA]-----------------------------
\appendix{A}
\label{sec:appa}
\setcounter{equation}{0}
\renewcommand{\theequation}{A.\arabic{equation}}
In many cases, if a kernel consists of a summation of components, 
the smoothness of the kernel with respect 
to a variable is dictated by a single component. One can then
define a so-called principal smoothness component: 
Let 
$\Theta\!=\!\{\vctr{y}-\vctr{x}\!\mid\!
\vctr{x}\!\in\!\underline\Theta,\vctr{y}\!\in\!\Theta\}$. 
If~$G^{\vctr{l}}(\vctr{t})=\sum_{{i}}G^{\vctr{l}}_{{i}}(\vctr{t})$ 
is asymptotically smooth in
${t}_k$ and there exist an index~${j}$, a~$\bar{\vctr{t}}\in\Theta$ and
a minimal order~$\bar{p}\in\IN$, 
such that for all~$\vctr{t}\in\Theta$, all
orders~$p\geq\bar{p}$ and all indices~$i$ it holds that:
\begin{equation}
\label{eq:prismoco}
\big|
G^{\vctr{l}-p\vctr{e}_k}_{i}(\vctr{t})
\big|
\leq
\alpha
\big|
G^{\vctr{l}-{p}\vctr{e}_k}_{j}(\bar{\vctr{t}}+(t_k-\bar{t}_k)\vctr{e}_k)
\big|
\:,
\end{equation}
for some positive constant~$\alpha$, then 
$\bar{G}{}^{\vctr{l}}_{k}(\eta)\equiv
\alpha\,{G}^{\vctr{l}}_{j}(\bar{\vctr{t}}+(\eta-\bar{t}_k)\vctr{e}_k)$ is 
the principal smoothness component of 
$G^{\vctr{l}}(\vctr{t})$ in the~$k$-direction.
Because the smoothness of the kernel in~$t_k$ is 
essentially governed by the 
principal smoothness component, suitable softening parameters,~$m_k$ and 
$p_k$, can be conveniently determined from the properties of 
$\bar{G}{}^{\vctr{l}}_{k}$. 

If the mesh width~$\vctr{h}$ is sufficiently small, then it is generally
possible to employ~$m_k=0$ and a fixed, minimal order of transfer, 
$\bar{p}_k$, depending only on~$l_k$, during several of the first 
coarsening stages. This is a result of the use of integrated kernels.
However, this does not apply if the evaluation is a
subtask in the solution of the integral equations, in which case
it is necessary to use~$m_k=O(\bar{p}_k)$ on the finest grids 
to maintain the stability properties of the single-grid operator.
For the isolated evaluation of the integral transform,  
\cite{Brandt:1998fk} shows that~$m_k=0$ requires~$\bar{p}_k>l_k$. Due to 
the singularity in the original kernel, the 
integrated kernel,~$G^{\vctr{l}}(\vctr{x},\vctr{y})$, 
contains components with singular 
or discontinuous derivatives at 
$y_k=x_k$ (\/$t_k=0$\/). Consequently, in order to determine the 
left-hand side of~\EQ{sreq}, it is necessary to distinguish between the 
region~$\Theta_s=\{\vctr{t}\in\Theta
\mid{}|t_k|\leq\bar{p}_k{}H_k/2\}$, {i.e.}
the region where  the singularity is in the interpolation interval, and the 
region~$\Theta\setminus\Theta_s$. In~$\Theta\setminus\Theta_s$, the 
kernel~$G^{\vctr{l}}(\vctr{x},\vctr{y})$ is~$\bar{p}_k$ times
differentiable with respect to~$x_k$ and~$y_k$. In this case the average order 
$\bar{p}_k$ interpolation error satisfies
\begin{equation}
\label{eq:interr}
\big\|
G^{\vctr{l}}_\vctr{H}- 
\II^{\vctr{H}-H_k\vctr{e}_k/2}_{\vctr{H}}G^{\vctr{\vctr{l}}}_\vctr{H}
\big\|_{1,\Theta\setminus\Theta_s}
\leq
\left( \gamma_2 H_k \right)^{\bar{p}_k}
\big\|G^{\vctr{l}(\bar{p}_k\vctr{e}_k)}\big\|_{1,\Theta\setminus\Theta_s}\:.
\end{equation}
Using~\EQ{prismoco},
\begin{equation}
\label{eq:GpscG}
\big\|G^{\vctr{l}(\bar{p}_k\vctr{e}_k)}\big\|_{1,\Theta\setminus\Theta_s}
\leq
\big\|\bar{G}^{{\vctr{l}}(\bar{p}_k)}_{k}\big\|_{1,\Theta\setminus\Theta_s}\:.
\end{equation}
Because~$\bar{G}{}^{\vctr{l}(\bar{p}_k-1)}_k$ 
vanishes asymptotically at infinity, the right-hand side of~\EQ{GpscG} is 
approximately~$|\bar{G}{}^{{\vctr{l}}(\bar{p}_k-1)}_{k}(\bar{p}_kH_k/2)|$. 
In~$\Theta_s$, however, equation~\EQ{interr} is useless as a result of the 
singular derivatives. The interpolation error in this region is bounded
by the the local value of the component with singular 
derivatives. Taylor expansion of the kernel in the neighborhood of the 
singularity yields that for sufficiently small~$t_k$ the behavior of this 
component is dominated by the principal smoothness component. 
Hence, the interpolation error for~$0\leq{}|y_k-x_k|\leq{}p_kH_k/2$ is 
bounded by
\begin{equation}
\label{eq:Ga=psca}
\big\|
G^{\vctr{l}}_\vctr{H}- 
\II^{\vctr{H}-H_k\vctr{e}_k/2}_{\vctr{H}}G^{\vctr{\vctr{l}}}_\vctr{H}
\big\|_{1,\Theta\setminus\Theta_s}
\leq
\big\|\bar{G}{}^{\vctr{l}}_{k}\big\|_{\text{max},\Theta_s}\:.
\end{equation}
Generally, the principal smoothness component is a monotonic function. The 
right-hand side of~\EQ{Ga=psca} is then
$|\bar{G}{}^{\vctr{l}}_{k}(\bar{p}_kH_k/2)|$. Notice that the contribution of 
the region~$\Theta_s$ to the entire integral is just~$O(\bar{p}_k H_k)$. 
By~\EQ{sreq} and~\EQ{interr} through~\EQ{Ga=psca}, 
it is anticipated that~$m_k=0$ can indeed be used on grids where the 
following two requirements are satisfied:
\begin{align}
\big|\bar{G}{}^{\vctr{l}}_{k}(\bar{p}_kH_k/2)\big|\,
\bar{p}_k H_k
&
\leq
h_k^{s_k-d(s_k-l_k)}
\,c_{G,u}\,c_h\,c_\gamma\:,
\label{eq:req1}\\
\big|\bar{G}{}^{\vctr{l}(\bar{p}_k-1)}_{k}(\bar{p}_kH_k/2)\big|\,
H_k^{\bar{p}_k}
&
\leq
h_k^{s_k-d(s_k-l_k)}
\,\gamma_2^{-\bar{p}_k}
\,c_{G,u}\,c_h\,c_\gamma\:.
\label{eq:req2}
\end{align}

Next, we investigate~$m_k>0$ and any~$p_k\geq\bar{p}_k$. The softened kernel 
is~$p_k-1$ times continuously differentiable in the~$k$-direction. Hence, the 
$p_k$-order interpolation error is bounded by
\begin{equation}
\label{eq:sinterr}
\big\|
G^{\vctr{l}}_\vctr{H}- 
\II^{\vctr{H}-H_k\vctr{e}_k/2}_{\vctr{H}}G^{\vctr{\vctr{l}}}_\vctr{H}
\big\|_{1,\Theta}
\leq
\left( \gamma_2 H_k \right)^{p_k}
\big\|G_\vctr{H}^{\vctr{l}(p_k\vctr{e}_k)}\big\|_{1,\Theta}\:.
\end{equation}
The derivative of the softened 
kernel in~\EQ{sinterr} is dominated by the derivative of its softened 
principal smoothness component:
\begin{equation}
\label{eq:sGpscG}
\big\|G_\vctr{H}^{\vctr{l}(p_k\vctr{e}_k)}\big\|_{1,\Theta}
\leq
\alpha_3\,
\big\|
\bar{G}_{\vctr{H}\,k}^{{\vctr{l}}(p_k)}\big\|_{1,\Theta}\:,
\end{equation}
for some positive constant~$\alpha_3$.
Due to the locality of the softening operation, it is necessary to distinguish 
two regions, viz., the 
region~$\Theta_s=\{\vctr{t}\in\Theta\mid |t_k|\leq{}(m_k+p_k/2)\,H_k\}$ 
where the softening domain is in the 
interpolation interval and the interpolation error is determined by the 
softening polynomial, and the 
region~$\Theta\setminus\Theta_s$, where the error results from the interpolation of 
the original kernel. In~$\Theta\setminus\Theta_s$, the average order~$p_k$ 
interpolation error is bounded by
\begin{equation}
\label{eq:sinterr2}
\big\|
G^{\vctr{l}}_\vctr{H}- 
\II^{\vctr{H}-H_k\vctr{e}_k/2}_{\vctr{H}}G^{\vctr{\vctr{l}}}_\vctr{H}
\big\|_{1,\Theta\setminus\Theta_s}
\leq
\left( \gamma_2 H_k \right)^{p_k}
\big\|\bar{G}{}_{k}^{{\vctr{l}}(p_k)}\big\|_{1,\Theta\setminus\Theta_s}\:.
\end{equation}
The right-hand side of~\EQ{sinterr2} is approximately 
$|\bar{G}{}^{{\vctr{l}}(p_k-1)}_{k}(m_kH_k)|$.
In the region~$\Theta_s$, the~$p_k$-derivative of the softening polynomial
determines the interpolation error. 
By~\EQ{soft1} and~\EQ{contcon}, for 
even (\/$odd=0$\/) and odd (\/$odd=1$\/) functions 
$\bar{G}{}_{k}^{{\vctr{l}}}$, the 
$p_k$-derivative of the softened principal smoothness component reads:
\begin{equation}
\label{eq:exPSC}
\bar{G}_{\vctr{H}\,k}^{{\vctr{l}}(p_k)}(\eta)=
\sum_{i=0}^{p_k-1}\sum_{j=0}^{p_k-1}
B_{ij}^{-1}
\,b(p_k,i)
\,(m_kH_k)^{-(2i+odd-j)}
\,\eta^{(2i+odd-p_k)}
\,\bar{G}{}^{\vctr{l}(j)}_k(m_kH_k)\:,
\end{equation}
for~$\eta\in[0,m_kH_k]$, with 
\begin{equation}
b(i,j)=\left\{
\begin{array}{ll}
{\displaystyle \frac{(2j+odd)!}{(2j+odd-i)!}}\qquad&2j+odd\geq{i}\:,\\
0&\mbox{otherwise}\:,
\end{array}
\right.
\end{equation}
and~$\mtrx{B}$ the~$(p_k\times{p_k})$-matrix with entries~$B_{ij}=b(i,j)$ 
(\/$i,j=0,\ldots,p_k-1$\/). The right-hand side of~\EQ{exPSC} can be used
to construct a convenient, sharp bound of the form:
\begin{equation}
\label{eq:shabo}
\alpha_3\,
\big\|
\bar{G}_{\vctr{H}\,k}^{{\vctr{l}}(p_k)}\big\|_{1,\Theta_s}
\leq
f(p_k,l_k)
\,F(m_kH_k)
\end{equation}
with~$F$ some elementary positive function.
The contribution of~$\Theta_s$ to the left-hand side of 
equation~\EQ{sreq} is just~$O(m_kH_k)$. By~\EQ{sreq} 
and~\EQ{sinterr} through~\EQ{shabo}, one arrives at the following two 
requirements for~$p_k$ and~$m_k$: 
\begin{align}
\left|\bar{G}{}^{\vctr{l}(p_k-1)}_{k}(m_kH_k)\right|
\left(H_k\right)^{p_k}
&
\leq
h_k^{s_k-d(s_k-l_k)}
\,c_{G,u}\,c_h
\,c_{\gamma}
\,\gamma_2^{-p_k}
\:,
\label{eq:req3}
\\
f(p_k,l_k)
\,m_kH_k^{p_k+1}
\,F(m_kH_k)
&
\leq
h_k^{s_k-d(s_k-l_k)}
\,c_{G,u}\,c_h
\,c_{\gamma}
\,\gamma_2^{-p_k}
\:.
\label{eq:req4}
\end{align}

Summarizing, whenever~$H_k$ satisfies~\EQ{req1} and~\EQ{req2},
$m_k=0$ and~$p_k=\bar{p}_k$ is used. Otherwise,~$m_k$ and~$p_k$ must
satisfy~\EQ{req3} and~\EQ{req4}.

%-----------------------[APPB]-----------------------------
\appendix{B}
\label{sec:appb}
\setcounter{equation}{0}
\renewcommand{\theequation}{B.\arabic{equation}}
To obtain the optimal transfer parameters for the fast evaluation of the 
model problem, the principal smoothness components of the 
kernel~\EQ{G22} is derived first. 
Because of the symmetry of~$G^{(2,2)}(\vctr{t})$ in 
${t}_1$ and~${t}_2$, it is sufficient to obtain the component in one 
direction. Observing that~$G^{(2,2)}_0(\vctr{t})$ in~\EQ{kerncom} 
can be recast into
\begin{equation}
\label{eq:T1psc}
G^{(2,2)}_0(\vctr{t})=
\frac{1}{2}{t}_{1}^2{t}_{2}\,\mbox{ln}\left({t}_2+\sqrt{{t}_1^2+{t}_2^2}\right)-
\frac{1}{2}{t}_{1}^2{t}_{2}\,\mbox{ln}\left|{t}_1\right|,
\end{equation}
analysis of the derivatives of the components reveals that the 
smoothness of~$G^{(2,2)}(\vctr{t})$ in the~$1$-direction is dominated by the 
second term in~\EQ{T1psc}. Assuming that the dimensions of the domains 
$\underline{\Omega}$ and~$\Omega$ are~$O(1)$, 
the principal smoothness component of~$G^{(2,2)}(\vctr{t})$ is 
\begin{equation}
\label{eq:psc22} 
\bar{G}{}_{k}^{(2,2)}(\eta)=\eta^2\mbox{ln}|\eta|.
\end{equation}
In general, the principal smoothness component in the~$k$-direction of 
$G^{\vctr{l}}(\vctr{t})$ for the family of kernels with 
$G(\vctr{x},\vctr{y})$ by~\EQ{modker} is
\begin{equation}
\label{eq:pscGa}
\bar{G}{}_{k}^{\vctr{l}}(\eta)=\eta^{l_k}\mbox{ln}|\eta|.
\end{equation}
Hence, it is useful to maintain a general notation. Because the principal 
smoothness component~\EQ{pscGa} is identical to 
the integrated kernel in the 1-dimensional model-problem 
treated in~\cite{Brandt:1998fk}, one can consult~\cite{Brandt:1998fk} for details of the
below optimization procedure.

By~\EQ{req1}, \EQ{req2} and~\EQ{pscGa},
it is anticipated that~$\bar{p}_k=l_k+2$ and~$m_{k}=0$, {i.e.} no softening 
at all, can be used on all grids with mesh size~$H_k$ satisfying
\begin{equation}
\label{eq:nosoft22}
H_k^{l_k+1}\mbox{ln}\left(H_{k}\right)=
O\left(h^{s_k-d(s_k-l_k)}\right).
\end{equation}

Next, consider~$m_{k}>0$. For all~$p_k\geq{}l_k+2$, requirement~\EQ{req4} 
is more restrictive than requirement~\EQ{req3} and, consequently,~$m_k$ and 
$p_k$ can be determined from the minimization of~\EQ{W_L} subject 
to~\EQ{req4}. From~\cite{Brandt:1998fk}, 
\begin{equation}
\big\|
\bar{G}_{\vctr{H}\,k}^{{\vctr{l}}(p_k)}\big\|_{1,\Theta_s}
\leq
f(p_k,l_k)
\,
(m_kH_k)^{l_k-p_k}\:,
\end{equation}
with~$f(p,l)=(2(p-l)/e)^{p-l}$.
Requirement~\EQ{req4} then assumes the following form:
\begin{equation}
\label{eq:reqspec}
(m_k/{\gamma_2})^{l_k-p_k+1}\,f(p_k,l_k)
\leq
g,
\end{equation}
where
\begin{equation}
g=
\frac{h^{s_k-d(s_k-l_k)}}{H_k^{l_k+1}}\,c_{G,u}\,c_h\,c_\gamma\,\gamma_2^{-(l_k+1)}\:.
\end{equation}
Assuming equality in~\EQ{reqspec}, 
\begin{equation}
\label{eq:mexp}
m_k=\gamma_2\,\mbox{exp}
\left( \frac{\mbox{ln}(g)-\mbox{ln}(f(p_k,l_k))}{l_k-p_k+1}\right)\:.
\end{equation}
Subsequently, \EQ{mexp} is substituted in~\EQ{W_L} and the~$p_k$ for 
which~$\mathrm{d}W^\vctr{H}/\mathrm{d}p_k=0$ is calculated. 
Making minor simplifications such as 
$(p_k-l_k)/(p_k-l_k-1)\approx{}1$, it follows that~$W^{\vctr{H}}$ is 
minimized when
\begin{equation}
\label{eq:e^x=x}
{\frac{32}{3}}\gamma_2\bar{H}^{-1} e^{1/\chi}=\chi,
\end{equation}
where
\begin{equation}
\label{eq:chi}
\chi=
\frac{-\left(p_k-l_k-1\right)}{\left(p_k-l_k-1\right)+\mbox{ln}(g)}.
\end{equation}

Summarizing, for the evaluation of the model transform, 
the optimal value of~$p_k$ for the transfer of the grid 
$\vctr{H}/2$ softened transform to grid~$\vctr{H}$ is the maximum of the 
lowest non-negative integer satisfying
\begin{equation}
\label{eq:pk}
p_k\geq\frac{-\chi}{\chi+1}\mbox{ln}(g)+l_k+1
\end{equation}
and~$p_k\geq{}l_k+2$. The corresponding optimal value of~$m_k$ is the 
first integer that satisfies:
\begin{equation}
\label{eq:mk}
m_k\geq\frac{3}{32}\bar{H}\gamma_2^{-1}\chi\left(p_k-l_k-1\right).
\end{equation}
Notice that~$\mbox{ln}(g)$ can be conveniently rewritten as
\begin{equation}
\label{eq:lng}
\mbox{ln}(g)=c^{}_a+(s_k-d(s_k-l_k))\,\mbox{ln}(h)-
(l_k+1)\,\mbox{ln}(H_k),
\end{equation}
where~$c_a$ is a constant depending on the geometry of the interpolation, 
the order of discretization, the average kernel value in~$\Theta$ and the 
derivatives of~$u(\vctr{y})$. Decreasing~$c_a$ 
causes~$p_k$ and~$m_k$ to increase, so that~$c_a$ controls the 
accuracy of the fast evaluation.

\end{document}